\documentclass{IOS-Book-Article}
\usepackage{amsmath}
\usepackage{amssymb}
\usepackage{amstext}
\usepackage{amsbsy}
\usepackage{amsopn}
\usepackage{amsthm}
\usepackage{amscd}
\usepackage{amsxtra}
\usepackage{amsfonts}
\usepackage{ifthen}
\usepackage{xspace}
\usepackage{fancyheadings}
\usepackage{layout}
\usepackage{hyperref}

\newcommand{\preprintserver}[2]{\href{http://arXiv.org/abs/math/#2}{#1/#2}}

\input xy
\xyoption{matrix}
\xyoption{arrow}
\xyoption{curve}
\newdir{ (}{{}*!/-5pt/@^{(}}
\newcommand{\xycenter}[1]{\begin{center}
                          \mbox{\xymatrix{#1}}
                          \end{center}
                         }

\newboolean{xlabels}
\newcommand{\xlabel}[1]{
                        \label{#1}
                        \ifthenelse{\boolean{xlabels}}
                                   {\marginpar[\hfill{\tiny #1}]{{\tiny #1}}}
                                   {}
                       }

\newcommand{\ZZ}{\mathbb{Z}}

\newcommand{\PZ}{\mathbb{P}}
\newcommand{\pz}{\mathbb{P}}

\newcommand{\AZ}{\mathbb{A}}

\newcommand{\CC}{\mathbb{C}}
\newcommand{\RR}{\mathbb{R}}
\newcommand{\QQ}{\mathbb{Q}}

\newcommand{\PP}{\mathbb{P}}
\newcommand{\GG}{\mathbb{G}}
\newcommand{\FF}{\mathbb{F}}

\newcommand{\sF}{{\mathcal F}}
\newcommand{\sG}{{\mathcal G}}

\newcommand{\sO}{{\mathcal O}}

\newcommand{\suchthat}{\, | \,}

\newboolean{probleme}
\newcommand{\problem}[1]
           {\ifthenelse{\boolean{probleme}}
                       {{\bf(PROBLEM: #1)\bf}}
                       {}
           }
\newboolean{zukuenftiges}
\newcommand{\zukunft}[1]
           {\ifthenelse{\boolean{zukuenftiges}}
                       {{\bf(AUSBAUM\"OGLICHKEIT: #1)\bf}}
                       {}
           }
\newboolean{extras}
\newcommand{\extra}[1]
           {\ifthenelse{\boolean{extras}}
                       {{\bf EXTRA #1 EXTRA\bf}}
                       {}
           }

\newboolean{ignore}
\newcommand{\ignore}[1]
           {\ifthenelse{\boolean{ignore}}
                       {{\bf IGNORE #1 IGNORE\bf}}
                       {}
           }

\DeclareMathOperator{\codim}{codim}

\DeclareMathOperator{\Hilb}{Hilb}

\DeclareMathOperator{\rank}{rank}

\DeclareMathOperator{\spec}{spec}

\DeclareMathOperator{\rad}{rad}

\theoremstyle{plain}
\begingroup
\newtheorem{thm}{Theorem}%[subsection]
\newtheorem{cor}[thm]{Corollary}
\newtheorem{lem}[thm]{Lemma}

\newtheorem{prop}[thm]{Proposition}

\numberwithin{thm}{subsection} 
\endgroup

\newtheorem*{thm*}{Theorem}
\newtheorem*{conj*}{Conjecture}
\newtheorem*{verm*}{Vermutung}

\theoremstyle{definition}
\newtheorem{defn}[thm]{Definition}
\newtheorem{rem}[thm]{Remark}
\newtheorem{example}[thm]{Example}

\newtheorem{experiment}[thm]{Experiment}
\newtheorem{alg}[thm]{Algorithm}
\newtheorem{heu}[thm]{Heuristic}

\numberwithin{equation}{section}

\newcommand{\nosubsections}{\renewcommand{\thethm}{\thesection.\arabic{thm}}
                            \setcounter{thm}{0}
                           }

\newcommand{\cref}[3]{(\ref{#1}, #2 \ref{#3})}

\date{\today}
\usepackage{amssymb,amsmath}
\setboolean{probleme}{true}
\setboolean{xlabels}{true}
\usepackage{graphicx}
\usepackage{verbatim}

\newcommand{\magmaExperiment}[2]{
	\newtheorem*{expB#1#2}{Experiment B.#1.#2}
	}

\magmaExperiment{1}{1}
\magmaExperiment{1}{2}
\magmaExperiment{1}{14}
\magmaExperiment{1}{18}
\magmaExperiment{2}{1}
\magmaExperiment{2}{4}
\magmaExperiment{2}{9}
\magmaExperiment{4}{3}
\magmaExperiment{4}{6}
\magmaExperiment{4}{10}
\magmaExperiment{5}{1}
\magmaExperiment{5}{2}
\magmaExperiment{5}{3}
\magmaExperiment{5}{4}

\begin{document}
\pagestyle{plain}
\begin{frontmatter}
\title{Finite Field Experiments}
\runningtitle{Finite Field Experiments}

\address{Institiut f\"ur Algebraische Geometrie\\ 
          Leibniz Universit\"at Hannover\\ 
          Welfengarten 1\\ 
          D-30167 Hannnover 
         }

%\email{\secemail}

%\urladdr{http://www-ifm.math.uni-hannover.de/\textasciitilde bothmer}

%\thanks{}

\author{\fnms{Hans-Christian} \snm{Graf v. Bothmer}}
\runningauthor{Hans-Christian Graf v. Bothmer}

\begin{abstract}
We show how to use experiments over finite fields to gain information about the solution set of polynomial equations in characteristic zero.
\end{abstract}
\end{frontmatter}

 %%%%%%%%%%%%%%%%%%%%%%%%%%%%%%%%%%%%%%%%%
 \section*{Introduction}
 %%%%%%%%%%%%%%%%%%%%%%%%%%%%%%%%%%%%%%%%%
 \nosubsections
 
 Let $X$ be a variety defined over $\ZZ$. According to Grothendieck we can picture $X$ as a family of varieties $X_p$ over $\spec \ZZ$ with fibers over closed points of $\spec \ZZ$ corresponding to reductions modulo $p$ and the generic fiber over $(0)$ corresponding to the variety $X_\QQ$ defined by the equations of $X$ over $\QQ$.

\begin{figure}[h]
\includegraphics*[width=10cm]{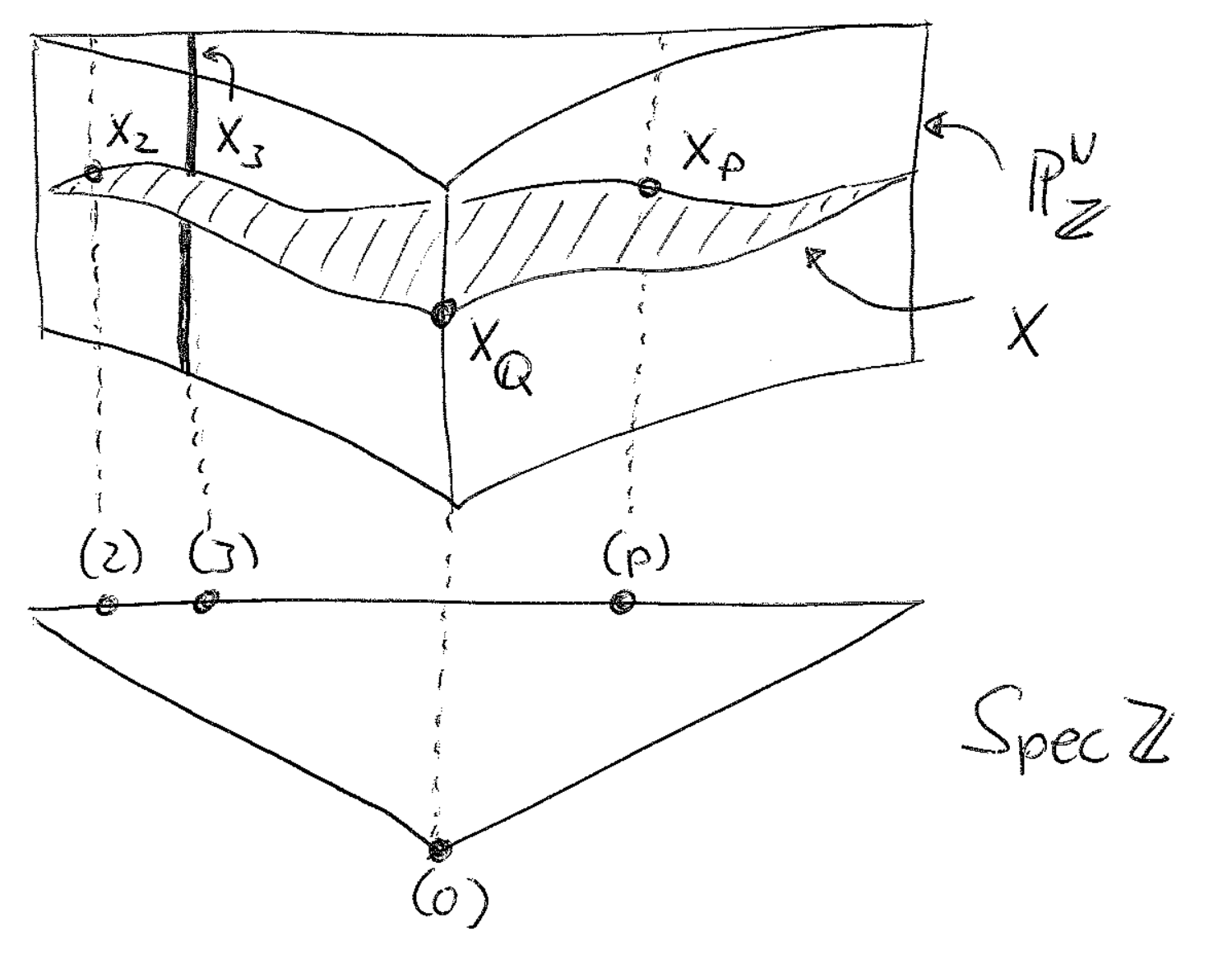}
\caption{A variety over $\spec \ZZ$}
\end{figure}

The generic fiber is related to the special fibers by semicontinuity theorems. For example, the dimension of $X_p$ is upper semicontinuous with 
$$\dim X_\QQ = \min_{p>0} \dim X_p.$$
This allows us to gain information about $X_\QQ$ by investigating $X_p$ which is often computationally much simpler.

Even more surprising is the relation between the geometry of $X_p$ and the number of $\FF_{p}$ rational points of $X_p$ discovered by Weil:

\begin{thm}
Let $X_p \subset \PP^n_{\FF_p}$ be a smooth curve of genus $g$, and $N$ be the number of $\FF_{p}$-rational points of $X_p$. Then
\[
		| 1-N + p | \le 2g \sqrt{p}.
\]
\end{thm}

He conjectured even more precise relations for varieties of arbitrary dimension which were proved by Deligne using $l$-adic cohomology.

In this tutorial we will use methods which are inspired by Weil's ideas, but are not nearly as deep. Rather we will rely on some basic probabilistic estimates which are nevertheless quite useful. I have learned these ideas from my advisor Frank Schreyer, but similar methods have been used independently by other people, for example Joachim von zur Gathen and Igor Shparlinski  \cite{ShaparlinskiComponents}, Oliver Labs \cite{labsPhD} and Noam Elkies \cite{ElkiesOberwolfach}.

The structure of these notes is as follows: We start in Section \ref{sGuessing} by evaluating 
the polynomials defining a variety $X$ at random points. This can give some heuristic information about the codimension $c$ of $X$ and about the number $d$ of codimension-$c$ components of $X$. 

In Section 
\ref{sTangents} we refine this method by looking at the tangent spaces of $X$ in random points. This
gives a way to also estimate the number of components in every codimension. As an application we
show how this can be applied to gain new information about the Poincar\'e center problem. 

In Section \ref{sExLift} we explain how it is often possible to prove that a solution
found over $\FF_p$ actually lifts to $\overline{\QQ}$. This is applied to the construction
of new surfaces in $\PP^4$. 

Often one would like not only to prove the existence of a lift, but explicitly find one. It is explained in Section \ref{sFindLift} how this can be
done if the solution set is zero dimensional. 

We close in Section \ref{sNodes} with a beautiful application of these lifting techniques
found by Oliver Labs, showing how he constructed a new septic with $99$ real nodes in $\PP^3_{\RR}$.

For all experiments in this tutorial we have used the computer algebra system Macaulay 2 \cite{M2}. The most important Macaulay 2 commands used are explained in Appendix \ref{Amacaulay}, for more detailed information we refer to the online help of Macaulay 2 \cite{M2}. In Appendix \ref{Amagma} Stefan Wiedmann provides a 
MAGMA translation of the Macualay 2 scripts in this tutorial. All scripts are available online at \cite{natoweb}. We would like to include translations to other computer algebra packages, so if you are for example a \verb|Singular|-expert, please contact us. 

Finally I would like to thank the referee for many valuable suggestions.

%%%%%%%%%%%%%%%%%%%%%%%%%%%%%%%%%%%%%%%%%%%
\section{Guessing} \label{sGuessing}
%%%%%%%%%%%%%%%%%%%%%%%%%%%%%%%%%%%%%%%%%%%
\nosubsections

\newcommand{\Fp}{\FF_p}
\newcommand{\Fpn}{\Fp[x_1,\dots,x_n]}
\newcommand{\term}[1]{{\sl #1}}
\renewcommand{\labelenumi}{(\roman{enumi})}
 
We start by considering the most simple case, namely that  of a hypersurface $X \subset \AZ^n$ defined by a
single polynomial $f \in \Fpn$. If $a \in \AZ^n$ is a point we have
\[
	f(a) = \left\{ 
	\begin{array}{cl}
	0 & \text{one possibility} \\
	\not = 0 & \text{$(p-1)$ possibilities} \\
	\end{array}
	\right.
\]
Naively we would therefore expect that we obtain zero for about $\frac{1}{p}$ of the points. 

\begin{experiment} \label{eSinglePolynomial}
We evaluate a given polynomial in $700$ random points, using Macaulay 2:
\begin{verbatim}
  R = ZZ[x,y,z,w]                  -- work in AA^4
  F = x^23+1248*y*z+w+129269698	   -- a Polynomial
  K = ZZ/7                         -- work over F_7
  L = apply(700,                   -- substitute 700 
       i->sub(F,random(K^1,K^4)))  -- random points
  tally L	    	      	             -- count the results
\end{verbatim}
obtaining:
\begin{verbatim}
  o5 = Tally{-1 => 100}
             -2 => 108
             -3 => 91
             0 => 98
             1 => 102
             2 => 101
             3 => 100
\end{verbatim}
Indeed, all elements of $\FF_7$ occur about $700/7 = 100$ times as one would expect naively. 
\end{experiment}

If $f = g\cdot h \in \Fpn$ is a reducible polynomial we have 
\[
	f(a) = g(a)h(a) = \left\{ 
	\begin{array}{cl}
	0 \cdot 0 & \text{$1$ possibility} \\
	* \cdot 0 & \text{$(p-1)$ possibilities} \\
	0  \cdot * & \text{$(p-1)$ possibilities} \\
	* \cdot * & \text{$(p-1)^2$ possibilities}
	\end{array}
	\right.
\]
so one might expect a zero for about $\frac{2p-1}{p^2} \approx \frac{2}{p}$ of the points. 

\begin{experiment} \label{eProduct}
We continue Experiment \ref{eSinglePolynomial} and
evaluate a product of two polynomials in $700$ random points:
\begin{verbatim}
  G = x*y*z*w+z^25-938493+x-z*w    -- a second polynomial
  tally apply(700,     	      	    -- substitute 700
     i->sub(F*G,random(K^1,K^4))) -- random points & count
\end{verbatim}
This gives:
\begin{verbatim}
  o8 = Tally{-1 => 86}
             -2 => 87
             -3 => 77
             0 => 198
             1 => 69
             2 => 84
             3 => 99
\end{verbatim}
Indeed, the value $0$ now occurs about twice as often, i.e. $198 \approx\frac{2}{7} \cdot 700$. 
\end{experiment}

\begin{figure}[h!] 
\includegraphics*[width=10cm]{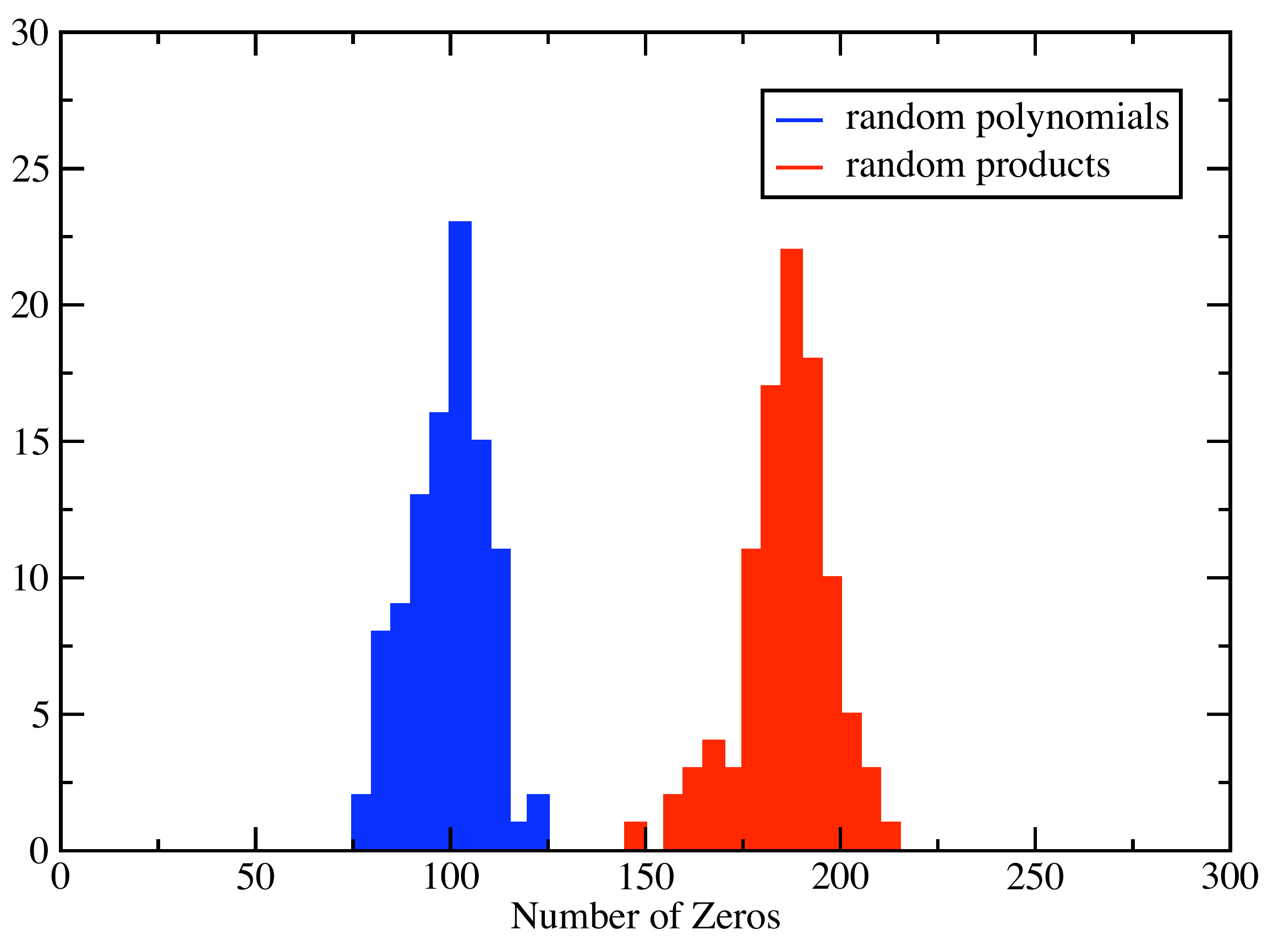}
\caption{Evaluating $100$ random polynomials and $100$ random products at $700$ points each.}
\label{fExperiment100}
\end{figure}

Repeating Experiments \ref{eSinglePolynomial} and \ref{eProduct} for $100$ random polynomials and
$100$ random products we obtain Figure \ref{fExperiment100}. 

Observe that the results for irreducible
and reducible polynomials do not overlap. Evaluating a polynomial at random points might therefore give some indication on the number of its irreducible factors. For this we will make the above naive observations more precise.

\begin{defn}
If $f \in \Fpn$ is a polynomial, we call the map
\[
	\begin{matrix}
	f|_{\Fp^n} &\colon& \Fp^n &\to &\Fp\\
	&& a &\mapsto&f(a)
	\end{matrix}
\]
the corresponding \term{polynomial function}. We denote by
\[
	V_p := \{ f \colon \Fp^n \to \Fp\}
\]
the vector space of all polynomial functions on $\Fp^n$.
\end{defn}

Being a polynomial function is nothing special:

\begin{lem} [Interpolation] \label{lInterpol}
Let $\phi \colon \Fp^n \to \Fp$ be any function. Then there exists a polynomial $f \in \Fpn$ such that
$\phi = f|_{\Fp^n}$.
\end{lem}

\begin{proof}
Notice that $(1-x^{p-1})=0 \iff x\not=0$. For every $a \in \Fp^n$ we define
\[
	f_a(x) := \prod_{i=1}^n (1-(x_i-a_i)^{p-1})
\]
and obtain
\[
	f_a(x) = \left\{ \begin{array}{cl}
		 1 & \text{if $x=a$} \\
		 0 & \text{if $x\not=a.$}
		 \end{array}
		 \right.
\]
Since $\Fp^n$ is finite we can consider $f := \sum_{a \in \Fp^n} \phi(a) f_a$ and obtain $f(x) = \phi(x)$ for
all $x \in \Fp^n$.
\end{proof}

\begin{rem} From Lemma \ref{lInterpol} it follows that
\begin{enumerate}
\item $V_p$ is a vector space of dimension $p^n$.
\item $V_p$ is a finite set with $p^{p^n}$ elements.
\item Two distinct polynomials can define the same polynomial function, for example $x^p$ and $x$. More generally if $F \colon \Fp \to \Fp$ is the \term{Frobenius endomorphism} then $f(a) = f(F(a))$ for
all polynomials $f$ and all $a \in \FF_p^n$. 
\end{enumerate}
\end{rem}

This makes it easy to count polynomial functions:

\begin{prop} \label{pDistributionSingle}
The number of polynomial functions $f \in V_p$ with $k$ zeros is
\[
	\binom{p^n}{k} \cdot 1^k \cdot (p-1)^{p^n-k}.
\]
\end{prop}

\begin{proof}
Since $V_p$ is simply the set of all functions $f \colon \Fp^n \to \Fp$, we can enumerate the ones with $k$  zeros as follows: First choose $k$ points and assign the value $0$ and then chose any of the other $(p-1)$ values for the remaining $p^n-k$ points.\
\end{proof}

\begin{cor}
The average number of zeros for polynomial functions $f \in V_p$ is
\[
	\mu= p^{n-1}
\]
and the standard deviation of the number of zeros in this set is 
\[
	\sigma = \sqrt{p^n\left(\frac{1}{p}\right)\left(\frac{p-1}{p}\right)}%= \sqrt{p^{n-1}-p^{n-2}}
	 < \sqrt{\mu}.
\]	
\end{cor}

\begin{proof}
Standard facts about binomial distributions.
\end{proof}

\begin{rem}
Using the normal approximation of the binomial distribution, we can estimate that more than $99\%$ of
all $f \in V_p$ satisfy
\[
	|\#V(f)-\mu|\le 2.58\sqrt{\mu}
\]
\end{rem}

For products of polynomials we have

\begin{prop} \label{pDistributionProduct}
The number of pairs  $(f,g) \in V_p \times V_p$ whose product has $k$ zeros is
\[
	\#\bigl\{(f,g) \bigr. \bigl| \#V(f\cdot g) = k\bigr\} = \binom{p^n}{k} \cdot (2p-1)^k \cdot ((p-1)^2)^{p^n-k}.
\]
In particular, the average number of zeros in this set is
\[
	\mu' = p^n\left(\frac{2p-1}{p^2}\right) %= 2q^{n-1} - q^{n-2} 
	\approx 2\mu
\]
and the standard deviation is
\[
	\sigma' =  \sqrt{p^n\left(\frac{2p-1}{p^2}\right)\left(\frac{(p-1)^2}{p^2}\right)} < \sqrt{\mu'}
\]

\begin{proof}
As in the proof of Proposition \ref{pDistributionSingle} we first choose $k$ points. For each of 
these points $x$ we choose either the value of $f(x) = 0$ and $g(x) \not=0$ or $f(x) \not=0$ and $g(x) =0$ or $f(x)=g(x)=0$. This gives $2p-1$ possibilities. For the remaining $p^n-k$ we choose $f$ and $q$ nonzero. For this we have $(p-1)^2$ possibilities. The formulas then follow again from standard facts about binomial distributions.
\end{proof}

\end{prop}
\begin{figure}
\includegraphics*[width=10cm]{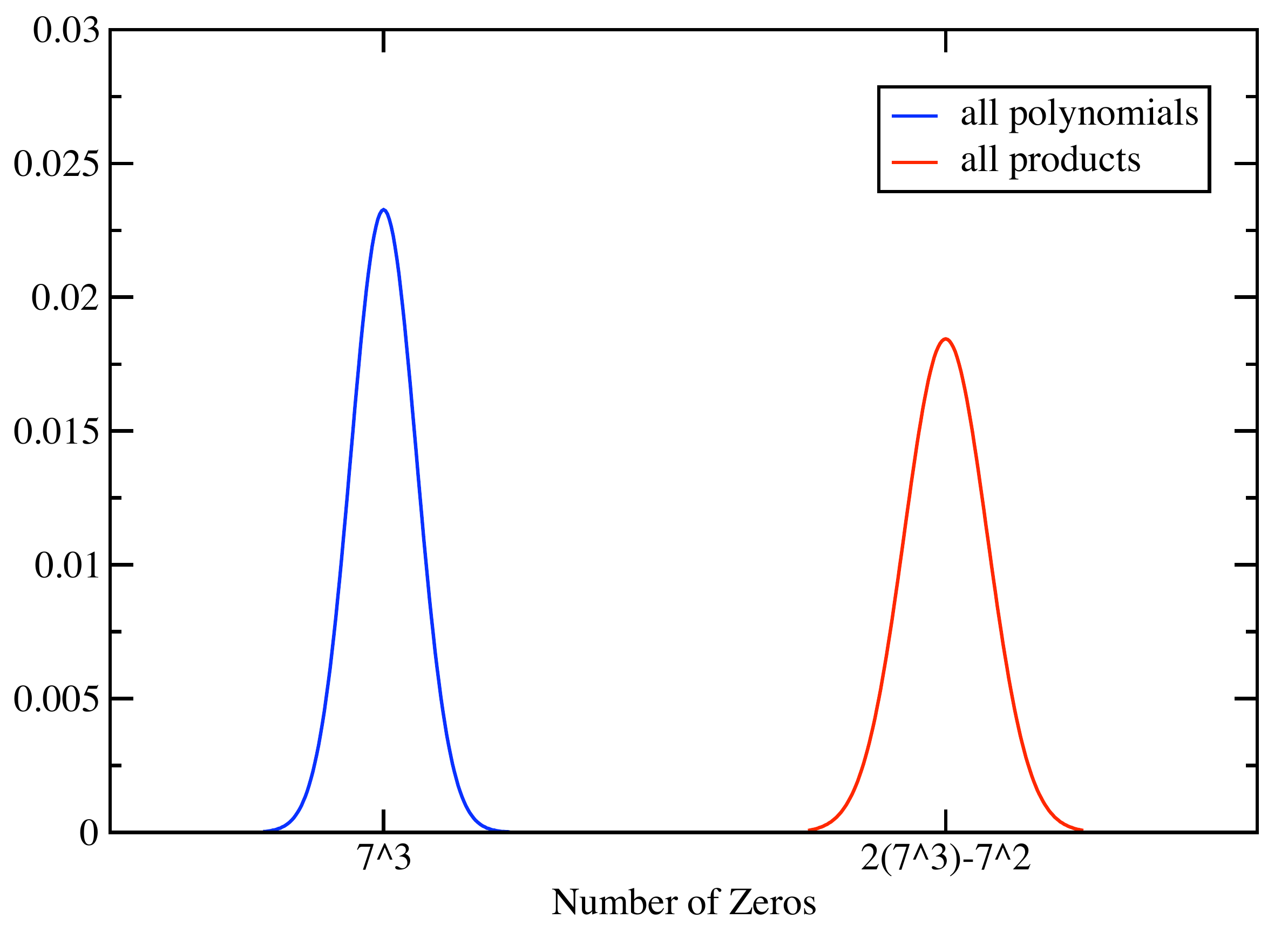}
\caption{Distribution of the number of zeros on hypersurfaces in $\AZ^4$ in characteristic $7$.}
\label{fCompare}
\end{figure}

\begin{rem}
It follows that more than $99\%$ of pairs $(f,g) \in V_p \times V_p$ 
satisfy
\[
	|\#V(f\cdot g)-\mu'|\le 2.58\sqrt{\mu'}.
\]
In particular, if  a polynomial $f$ has a number of zeros that lies outside of this range one can reject the hypothesis that $f$ is a product of two irreducible with $99\%$ confidence.
\end{rem}

Even for small $p$ the distributions of Proposition \ref{pDistributionSingle} and Proposition \ref{pDistributionProduct} differ substantially (see Figure \ref{fCompare}).

\begin{figure}
\includegraphics*[width=10cm]{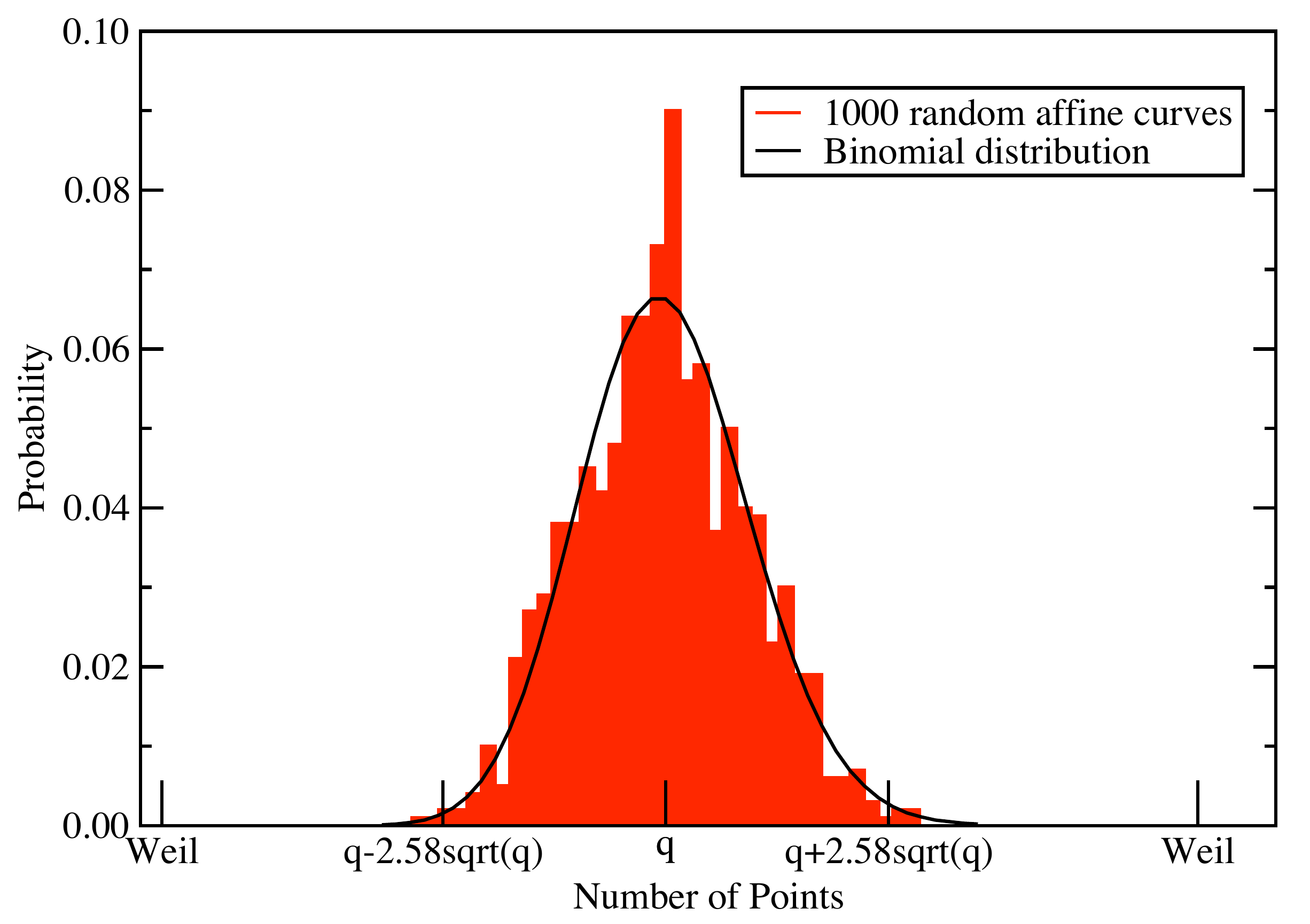}
\caption{Number of points on $1000$ affine quartics in characteristic $q=37$ compared to
the corresponding binomial distribution and Weil's bound.}
\label{fPlaneCurves}
\end{figure}

\begin{rem}
For plane curves we can compare our result to the Weil conjectures. Weil shows that $100\%$ of smooth pane curves of genus $g$ in $\PP^2_{\FF_p}$ satisfy
\[
	| N-(p+1) | \le 2g \sqrt{p}
\]
while we proved that $99\%$ of the polynomial functions on $\AZ^2$ satisfy
\[
	| N - p | \le 2.58 \sqrt{p}.
\]
Of course Weil's theorem is much stronger. If $p > 4g^2$ Weil's theorem implies for example that
every smooth curve of genus $g$ over $\FF_p$ has a rational point, while no such statement can be derived from our results. If on the other hand one is satisfied with approximate results, our estimates have the advantage that they are independent of the genus $g$. In Figure \ref{fPlaneCurves} we compare the two results  with an experiment in the case of plane quartics.  (Notice that smooth plane quartics have genus $3$.)
\end{rem}

For big $n$ it is very time consuming to count all $\FF_p$-rational points on $V(f) \subset \AZ^n$. We can avoid this problem by using a statistical approach once again.

\newcommand{\gammahat}{\hat{\gamma}}

\begin{defn}
Let $X \subset \AZ^n$ be a variety over $\FF_p$. Then
\[
	\gamma_p(X) := \frac{\# X(\FF_p)}{\# \AZ^n(\FF_p)}
\]
is called the {\sl fraction of zeros} of $X$. If furthermore $x_1,\dots,x_m \in \AZ^n(\FF_p)$ are points
then
\[
	\gammahat_p(X) := \frac{\#\{ i \suchthat x_i \in X\}}{m}
\]
is called an {\sl empirical fraction of zeros} of $X$. 
\end{defn}

\begin{rem}
If we choose the points $x_i$ randomly and independently, the probability that $x_i \in X$ is $\gamma_p(X)$. Therefore
we have the following:
\begin{enumerate} 
\item $\mu(\gammahat_p) = \gamma_p$, i.e. for large $m$ we expect $\gammahat_p(X) \approx \gamma_p(X)$
\item $\sigma(\gammahat_p) \approx \sqrt{\frac{\gamma_p}{m}}$, i.e. the quadratic mean of the error $|\gamma_p-\gammahat_p|$ decreases with $\sqrt{m}$.
\item Since for hypersurfaces $\gamma_p \approx \frac{1}{p}$ the average error depends neither on the number of variables $n$ nor on the degree of $X$.
\item Using the normal approximation again one can show that it is usually enough to 
test about $100\cdot p$ points to distinguish between reducible and irreducible polynomials (for more precise estimates see \cite{irred}).
\end{enumerate}
\end{rem}

\begin{experiment} \label{eQuadrics1}
Consider quadrics in $\PP^3$ and let
$$
	\Delta := \{ \text{singular quadric} \} \subset \{ \text{all quadrics}\} \cong \PP^9
$$
be the subvariety of singular quadrics in the space of all quadrics. Since having a singularity is a codimension $1$ condition for surfaces in $\PP^3$ we expect $\Delta$ to be a hypersurface. Is $\Delta$ irreducible? 
Using our methods we obtain a heuristic answer using Macaulay 2:
\begin{verbatim}
  -- work in characteristic 7 
  K = ZZ/7
  -- the coordinate Ring of IP^3 
  R = K[x,y,z,w]
  -- look at 700 quadrics 
  tally apply(700, i->codim singularLocus(ideal random(2,R)))
\end{verbatim}
giving
\begin{verbatim}
  o12 = Tally{2 => 5  }
              3 => 89
              4 => 606.
\end{verbatim}
We see $95 = 89+5$ of our $700$ quadrics were singular, i.e. $\gammahat(\Delta) = \frac{95}{700}$.
Since this is much closer to $\frac{1}{7}$, then it is to $\frac{2}{7}$ we guess that $\Delta$ is irreducible. 
Notice that we have not even used the equation of $\Delta$ to obtain this estimate.
 \end{experiment}

\begin{figure}
\includegraphics*[width=6cm]{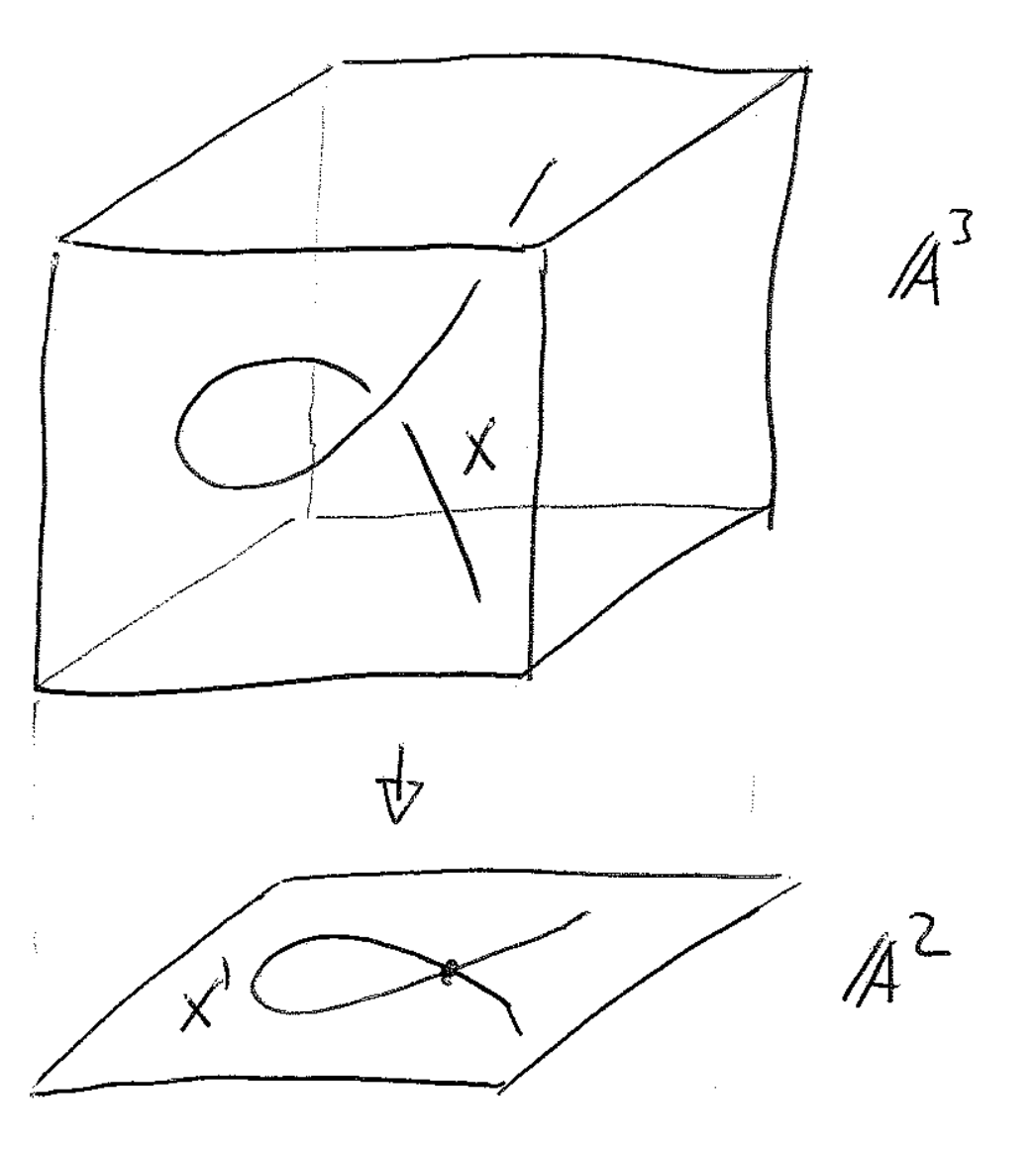}
\caption{The projection of a curve in $\AZ^3$ is a hypersurface in $\AZ^2$ and most
points have only one preimage.}
\label{fProjection}
\end{figure}

Let's now consider an irreducible variety $X \subset \AZ^n$  of codimension $c>1$. Projecting $\AZ^n$ to a subspace $\AZ^{n-c+1}$ we obtain a projection $X' \subset \AZ^{n-c+1}$ of $X$ (see Figure \ref{fProjection}). Generically $X'$ is a hypersurface, so by our arguments above $X'$ has approximately $p^{n-c}$  points. Generically most points of $X'$ have only one preimage in $X$ so we obtain the following very rough heuristic: 

\begin{heu} \label{hHigherCodim}
Let $X \subset \AZ^n_{\FF_p}$ be a variety of codimension $c$ and $d$ the number of components of codimension $c$, then 
\[
	\gammahat_p(X) \approx \frac{d}{p^c}
\]
\end{heu}

\begin{rem}
A more precise argument for this heuristic comes from the Weil Conjectures. Indeed, the
number of $\FF_p$-rational points on an absolutely irreducible projective variety $X$ is
\[
	p^{\dim X} + \text{lower order terms},
\]
so $\gamma_p \approx \frac{1}{p^{\codim X}}$. Our elementary arguments still work in the case of complete intersections and determinantal varieties \cite{irred}. 
\end{rem}

\begin{rem}
Notice that Heuristic \ref{hHigherCodim} involves two unknowns: $c$ and $d$. To determine these
one has to measure over several primes of good reduction.
\end{rem}

\begin{experiment}
As in Experiment \ref{eQuadrics1} we look at quadrics in $\PP^3$. These are given by their $10$ coefficients and form a $\PP^9$.  This time we are interested in the variety $X \subset \PP^9$ of quadrics whose singular locus is at least one dimensional. For this we first
define a function that looks at random quadrics over $\FF_p$ until it has found at least $k$ examples
whose singular locus has codimension at most $c$. It then returns the number of trials needed to do this.

\begin{verbatim}
  findk = (p,k,c) -> (
       K := ZZ/p;
       R := K[x,y,z,w];
       trials := 0;
       found := 0;
       while found < k do (
         Q := ideal random(2,R);
         if c>=codim (Q+ideal jacobian Q) then (
	           found = found + 1;
	           print found;
	           );
  	      trials = trials + 1;
  	      );
       trials
       )	       
\end{verbatim}

Here we use \verb!(Q+ideal jacobian Q)! instead of  \verb!singularLocus(Q)!, since the second option quickly produces a memory overflow. 

The function \verb!findk! is useful since the error in estimating $\gamma$ from $\gammahat$ depends on the number of singular quadrics found. By searching until a given number of singular quadrics is found make sure that the error estimates will be small enough. 

We now look for quadrics that have singularities of dimension at least one
\begin{verbatim}
  k=50; time L1 = apply({5,7,11},q->(q,time findk(q,k,2)))
\end{verbatim}
obtaining
\begin{verbatim}
  {(5, 5724), (7, 17825), (11, 68349)}
\end{verbatim}
i.e. $\gamma_5 \approx \frac{50}{5724}$, $\gamma_7 \approx \frac{50}{17825}$ and $\gamma_{11} \approx \frac{50}{68349}$. The codimension $c$ of $X$ can be interpreted as the negative slope
in a log-log plot of $\gamma_p(X)$ since Heuristic \ref{hHigherCodim} gives 
$$
\gammahat_p(X) \approx \frac{d}{p^c} \iff
	\log(\gammahat_p(X)) \approx \log(d) - c\log(p).
$$
This is illustrated in Figure \ref{fLogLog}. 

\begin{figure}
\includegraphics*[width=10cm]{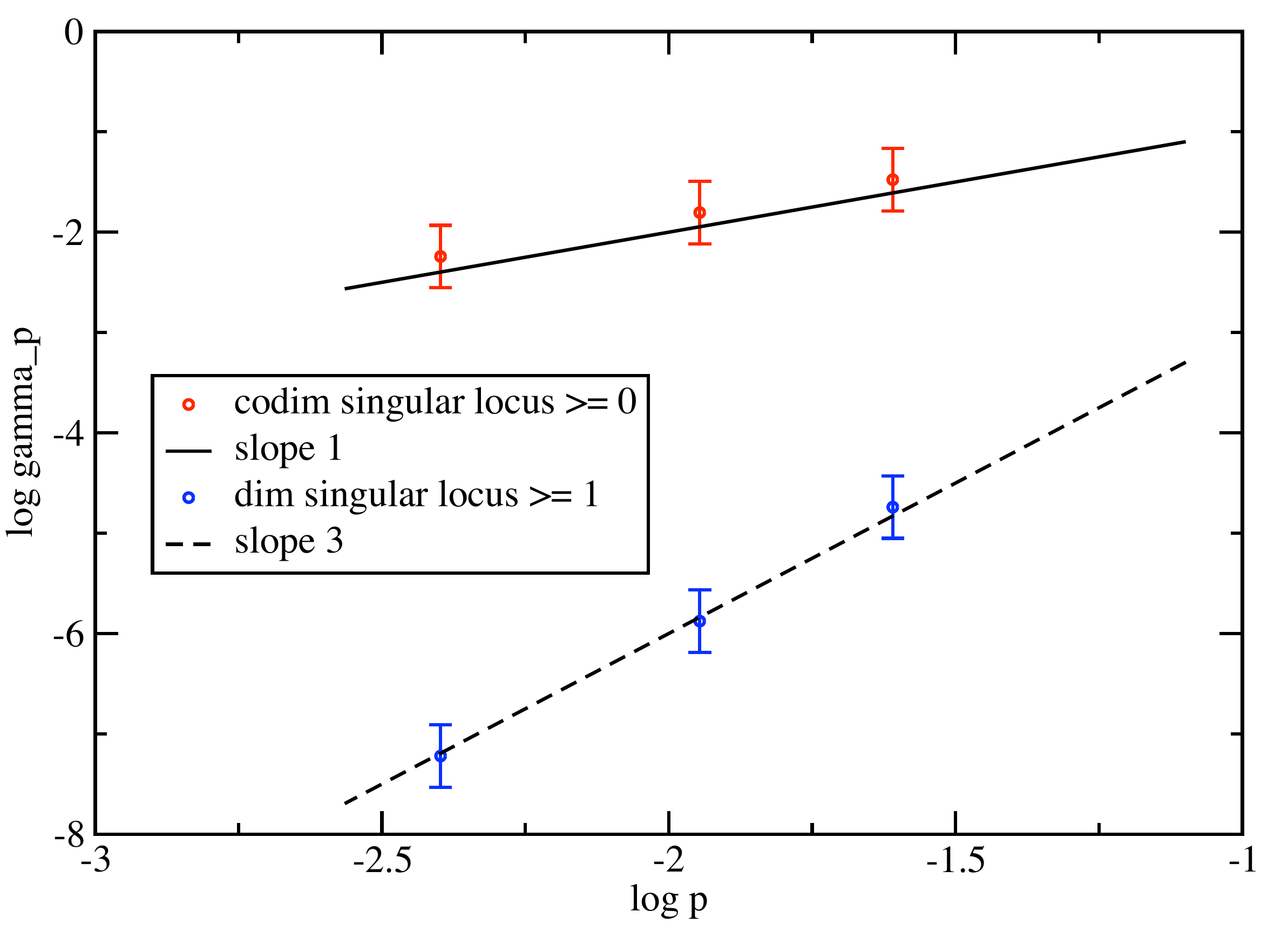}
\caption{Measuring the codimension of quadrics with zero and one dimensional singular loci. Here we compare the measurements with lines of the correct slope $1$ and $3$. }
\label{fLogLog}
\end{figure}

By using \verb!findk! with $k=50$ the errors of all our measurements 
are of the same magnitude.  We can
therefore use regression to calculate the slope of a line fitting these measurements:
\nocite{Bronstein}

\begin{verbatim}
  -- calculate slope of regression line by 
  -- formula from [2] p. 800 
  slope = (L) -> (
       xbar := sum(apply(L,l->l#0))/#L;
       ybar := sum(apply(L,l->l#1))/#L;
       sum(apply(L,l->(l#0-xbar)*(l#1-ybar)))/
         sum(apply(L,l->(l#0-xbar)^2))
       )

  -- slope for dim 1 singularities
  slope(apply(L1,l->(log(1/l#0),log(k/l#1))))

  o5 = 3.13578
\end{verbatim}
The codimension of $X$ is indeed $3$ as can be seen
by the following geometric argument: Each quadric with a singular locus of dimension $1$ is a union of two hyperplanes. Since the family $\hat{\PP^3}$ of all hyperplanes in $\PP^3$ is $3$-dimensional, we obtain $\dim X = 6$ which has codimension $3$ in the $\PP^9$ of all quadrics.
\end{experiment}

The approach presented in this section measures the number of components of minimal codimension quite well. At the same time it is very difficult to see components of larger codimension. One reason is 
%that the precision of $\gammahat$ improves only with $\sqrt{m}$. Another reason is 
that the rough approximations that we have made introduce errors in the order of $\frac{1}{p^{c+1}}$. 

We will see in the next section how one can circumvent these problems.

%EXCERSISE: Rank filtration of symmetric, antisymmetric and cataclectic matrices.
 
 %%%%%%%%%%%%%%%%%%%%%%%%%%%%%%%%%%%%%%%%%%%%%
 \section{Using Tangent Spaces} \label{sTangents}
 %%%%%%%%%%%%%%%%%%%%%%%%%%%%%%%%%%%%%%%%%%%%%
 \nosubsections
 
 If $X \subset \AZ^n$ has components of different dimensions, the guessing method of Section \ref{sGuessing} does not detect the smaller components.
 
\begin{figure}[h!]
\includegraphics*[width=6cm]{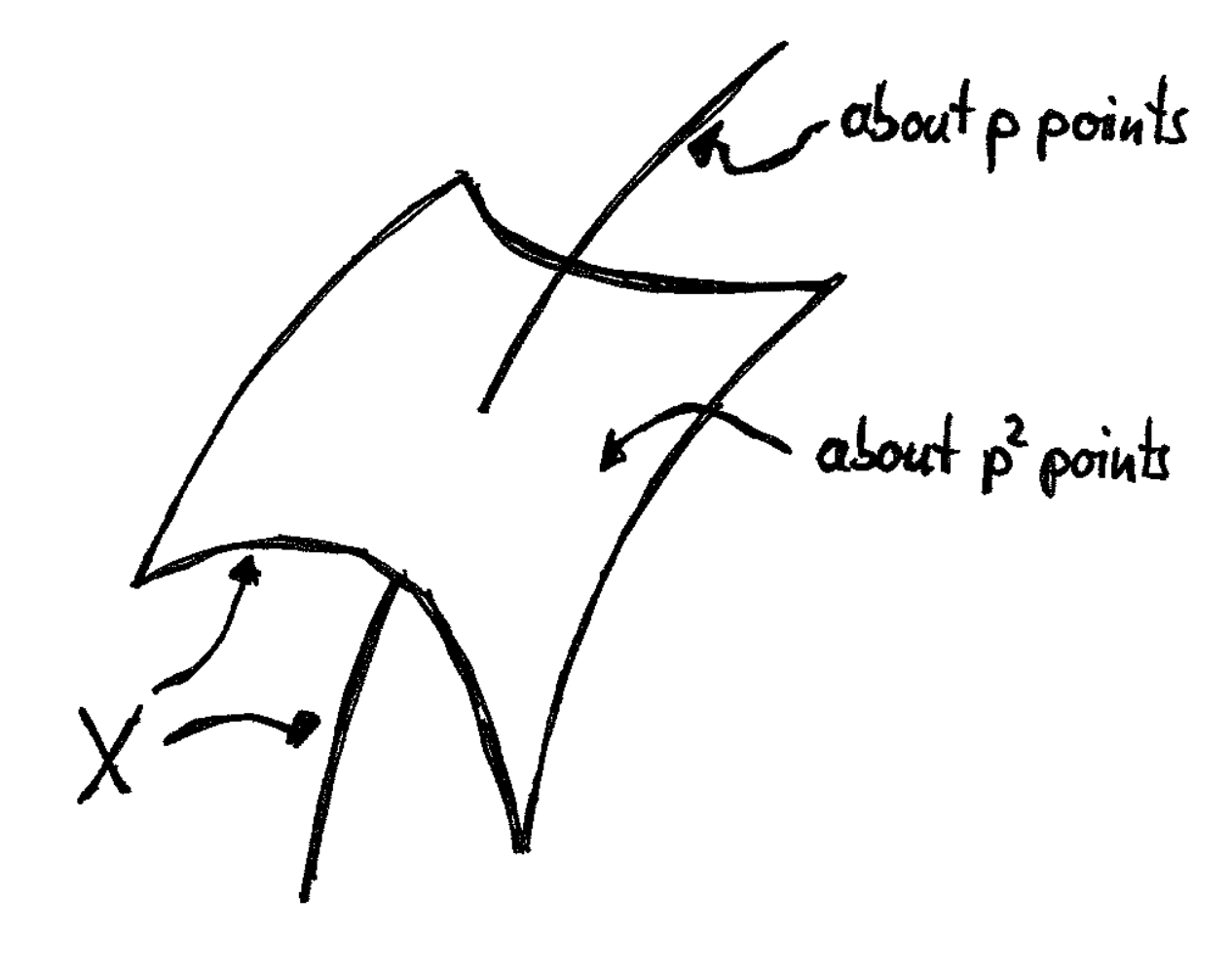}
\caption{Expected number of $\FF_p$ rational points on a union of a curve and a surface.}
\label{fDifferentDimPoints}
\end{figure}

 If  for example $X$ is the union of a curve and a surface in $\AZ^3$, we expect the surface to have about $p^2$ points while the curve will have about $p$ points (see Figure \ref{fDifferentDimPoints}). 
 
 Using Heuristic \ref{hHigherCodim} we obtain
 \[
 	\gamma_p = \frac{p^2+p}{p^3} \approx \frac{1}{p}
\]
indicating that $X$ has $1$ component of codimension $1$. The codimension $2$ component remains invisible.
 
\begin{experiment} \label{eFGH}
Let's check the above reasoning in an experiment. First define a function that produces
a random inhomogeneous polynomial of given degree:
\begin{verbatim}
  randomAffine = (d,R) -> sum apply(d+1,i->random(i,R))
\end{verbatim}
with this we choose random polynomials $F$, $G$ and $H$ in $6$ variables
\begin{verbatim}
  n=6
  R=ZZ[x_1..x_n];
  F = randomAffine(2,R)
  G = randomAffine(6,R);
  H = randomAffine(7,R);
\end{verbatim}
and consider the ideal $I = (FG,FH)$ % = (F) \cap (G,H)$
\begin{verbatim}
  I = ideal(F*G,F*H);
\end{verbatim}
Finally, we evaluate the polynomials of $I$ in $700$ points of characteristic $7$ and count
how many of them lie in $X = V(I)$:
\begin{verbatim}
  K = ZZ/7
  t = tally apply(700,i->(
  	  0 == sub(I,random(K^1,K^n))
  	  ))
\end{verbatim}
This yields
\begin{verbatim}
  o9 = Tally{false => 598}
             true => 102
\end{verbatim}
i.e. $\gammahat_7(X) = \frac{102}{700}$ which is very close to $\frac{1}{7}$. Consequently we would
conclude that $X$ has one component of codimension $1$. The codimension $2$ component  given by $G=H=0$ remains invisible. 
\end{experiment}

To improve this situation we will look at tangent spaces. Let $a \in X \subset \AZ^n$ be a point and $T_{X,a}$ the tangent space of $X$ in $a$. If $I_X = (f_1,\dots,f_m)$, let
 \[
 	J_X = \begin{pmatrix}
			\frac{d f_1}{d x_1} & \dots & \frac{d f_1}{d x_n}  \\
			\vdots & & \vdots \\
			\frac{d f_m}{d x_1} & \dots & \frac{d f_m}{d x_n}  \\
		\end{pmatrix}
\]
be the Jacobian matrix. We know from differential geometry that
 \[
 	T_{X,a} = \ker J_X(a) = \{v \in K^n \suchthat J_X(a)v = 0 \}.
\]
We can use tangent spaces to estimate the dimension of components of $X$: 

\begin{prop}
Let $a \in X \subset \AZ^n$ be a point and $X' \subset X$ a component containing $a$. Then
$\dim X' \le \dim T_{X,a}$ with equality holding in smooth points of $X$.
\end{prop}

\begin{proof} \cite[II.1.4. Theorem 3]{Shaf1} \end{proof}

\begin{figure}[h!]
 \includegraphics*[width=6cm]{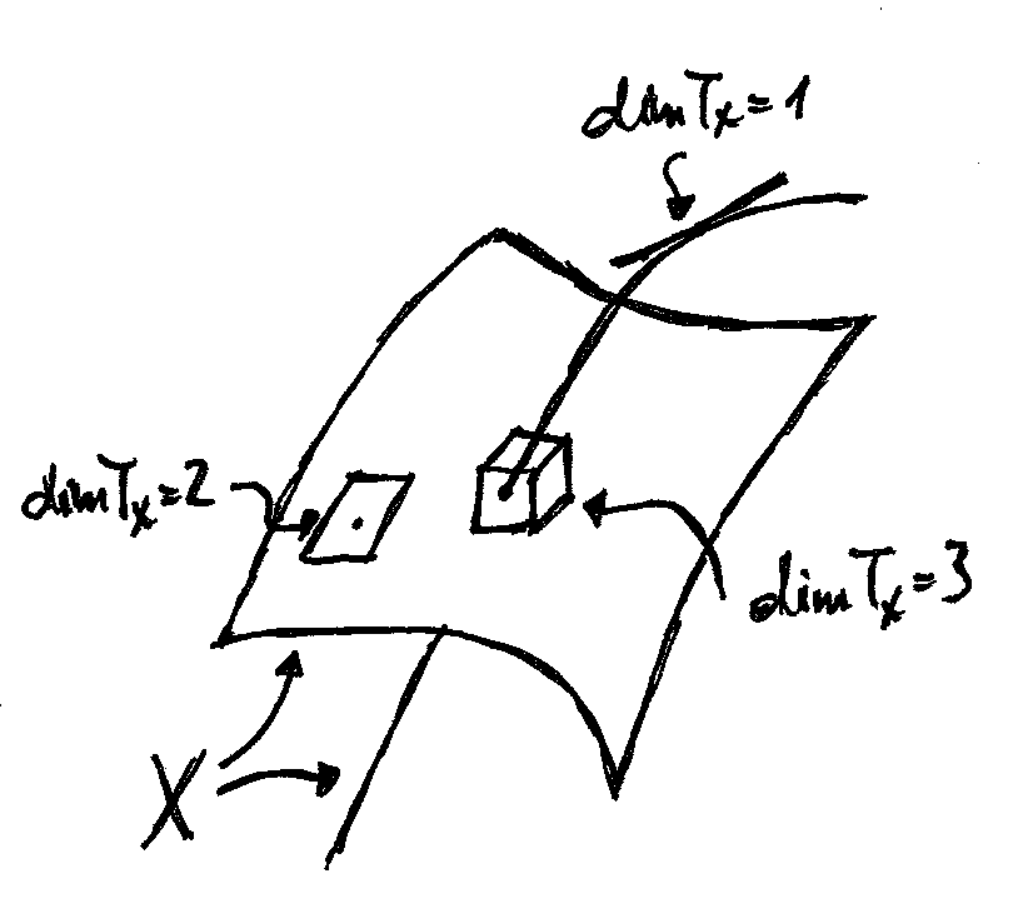}
 \caption{Dimension of tangent spaces in $\FF_p$ rational points on a union of a curve and a surface.}
 \label{fDifferentDimTangents}
 \end{figure}

In particular, we can use the dimension of the tangent space in a point $a \in X$ to
separate points that lie on different dimensional components, at least if these components are non reduced (see Figure \ref{fDifferentDimTangents}). For each of these sets we use Heuristic \ref{hHigherCodim} to obtain

\begin{heu} \label{hTangents}
Let $X \subset \AZ^n$ be a variety . If $J_X$ is the Jacobian matrix of $X$ and  $a_1,\dots,a_m\in \AZ^n$ are points, then the number of codimension $c$ components of $X$
is approximately 
\[
	\frac{\#\{i \suchthat \text{$a_i \in X$ and $\rank J_X(a_i) = c$}\}\cdot p^c}{m}
\]
\end{heu}

\begin{experiment} \label{eFGHtangent}
Let's test this heuristic by continuing Experiment \ref{eFGH}. For this we first calculate the
Jacobian matrix of the ideal $I$
\begin{verbatim}
  J = jacobian I;
\end{verbatim}
Now we check again $700$ random points, but when we find a point on $X=V(I)$ we also calculate the rank of the Jacobian matrix in this point:
\begin{verbatim}
  K=ZZ/7
  time t = tally apply(700,i->(
  	  point := random(K^1,K^n);
  	  if sub(I,point) == 0 then 
  	       rank sub(J,point)
  	  ))
\end{verbatim}
The result is
\begin{verbatim}
  o12 = Tally{0 => 2     }
              1 => 106
              2 => 14
              null => 578
\end{verbatim}
Indeed, we find that there are about $\frac{106 \cdot 7^1}{700} = 1.06$ components of dimension $1$ and about $\frac{14 \cdot 7^2}{700} = 0.98$ components of codimension $2$. For codimension $0$ the result is $\frac{2 \cdot 7^0}{700} \approx 0.003$ consistent with the fact that there are no components of codimension $0$.
\end{experiment}

\begin{rem}
It is a little dangerous to give the measurements as in Experiment \ref{eFGHtangent} without error bounds. Using the Poisson approximation of binomial distributions with small success probability we obtain 
$$\sigma(\text{number of points found}) \approx \sqrt{\text{number of points found}}.$$
In the above experiment this gives
\[
	\text{\# codim $1$ components} = 
	\frac{(106 \pm 2.58\sqrt{106}) \cdot 7^1}{700} = 1.06 \pm 0.27
\]
and
\[
	\text{\# codim $2$ components} = 
	\frac{(14 \pm 2.58\sqrt{14}) \cdot 7^2}{700} = 0.98 \pm 0.68.
\]
where the error terms denote the $99\%$ confidence interval.
Notice that the measurement of the codimension $2$ components is less precise. As a rule of thumb good error bounds are obtained if one searches until about $50$ to $100$ points of interest are found. 
\end{rem}

\begin{rem}
This heuristic assumes that the components do not intersect. If components do have high dimensional intersections, the heuristic might give too few components, since intersection points are singular and have lower codimensional tangent spaces. 
\end{rem}

In more involved examples calculating and storing the Jacobian matrix $J_X$ can use a lot of time and space. Fortunately one can  calculate $J_X(a)$ directly without calculating $J_X$ first:

\begin{prop}
Let $f \in \Fpn$ be a polynomial, $a\in \Fp^n$ a point and $b \in \Fp^n$ a vector. Then
\[
	f(a+ b \varepsilon) = f(a) + d_b f(a) \varepsilon \in \Fp^n[\varepsilon]/(\varepsilon^2).
\]
with $d_b f$ denoting the derivative of $f$ in direction of $b$. 
In particular, if $e_i \in \Fp^n$ is the $i$-th unit vector, we have 
$$f(a+e_i \varepsilon) = f(a) + \frac{d f}{d x_i}(a)\varepsilon.$$
\end{prop} 

\begin{proof} Use the Taylor expansion.
\end{proof}

\begin{example}
$f(x) = x^2 \implies f(1+\varepsilon) = (1+\varepsilon)^2 = 1+2\varepsilon = f(1) + \varepsilon f'(1)$
\end{example}

\begin{experiment}
To compare the two methods of calculating derivatives, we consider
the determinant of a random matrix with polynomial entries. First we create a random matrix
\begin{verbatim}
  K = ZZ/7                      -- characteristic 7
  R = K[x_1..x_6]               -- 6 variables
  M = random(R^{5:0},R^{5:-2})  -- a random 5x5 matrix with 
                                -- quadratic entries
\end{verbatim}
calculate the determinant
\begin{verbatim}
  time F = det M;
  -- used 13.3 seconds
\end{verbatim}
and its derivative with respect to $x_1$.
\begin{verbatim}
  time F1 = diff(x_1,F);
  -- used 0.01 seconds
\end{verbatim}
Now we substitute a random point:
\begin{verbatim}
  point = random(K^1,K^6)
  time sub(F1,point) 
  -- used 0. seconds
  
  o7 = 2
\end{verbatim}
By far the most time is used to calculate the determinant. With the $\varepsilon$-method this
can be avoided. We start by creating a vector in the direction of $x_1$:
\begin{verbatim}
  T = K[e]/(e^2)                       -- a ring with e^2=0
  e1 = matrix{{1,0,0,0,0,0}}           -- the first unit vector
  point1 = sub(point,T) + e*sub(e1,T)  -- point with direction
\end{verbatim}
Now we first evaluate the matrix $M$ in this vector
\begin{verbatim}
  time M1 = sub(M,point1)
  -- used 0. seconds
\end{verbatim}
and only then take the determinant
\begin{verbatim}
  time det sub(M,point1)
  -- used 0. seconds
  
  o12 = 2e + 1
\end{verbatim}
Indeed, the coefficient of $e$ is the derivative of the determinant in this point. This method is too fast to measure by the \verb=time= command of Macaulay 2. To get a better time estimate, we calculate the derivative of the determinant at 5000 random points:
\begin{verbatim}
  time apply(5000,i->(               
	  point := random(K^1,K^6);    -- random point
       	  point1 := sub(point,T)+e*sub(e1,T); -- tangent direction
	  det sub(M,point1);               -- calculate derivative
	  ));
  -- used 12.76 seconds
\end{verbatim}
Notice that this is still faster than calculating the complete determinant once.
\end{experiment}

\begin{figure} 
\includegraphics*[width=10cm]{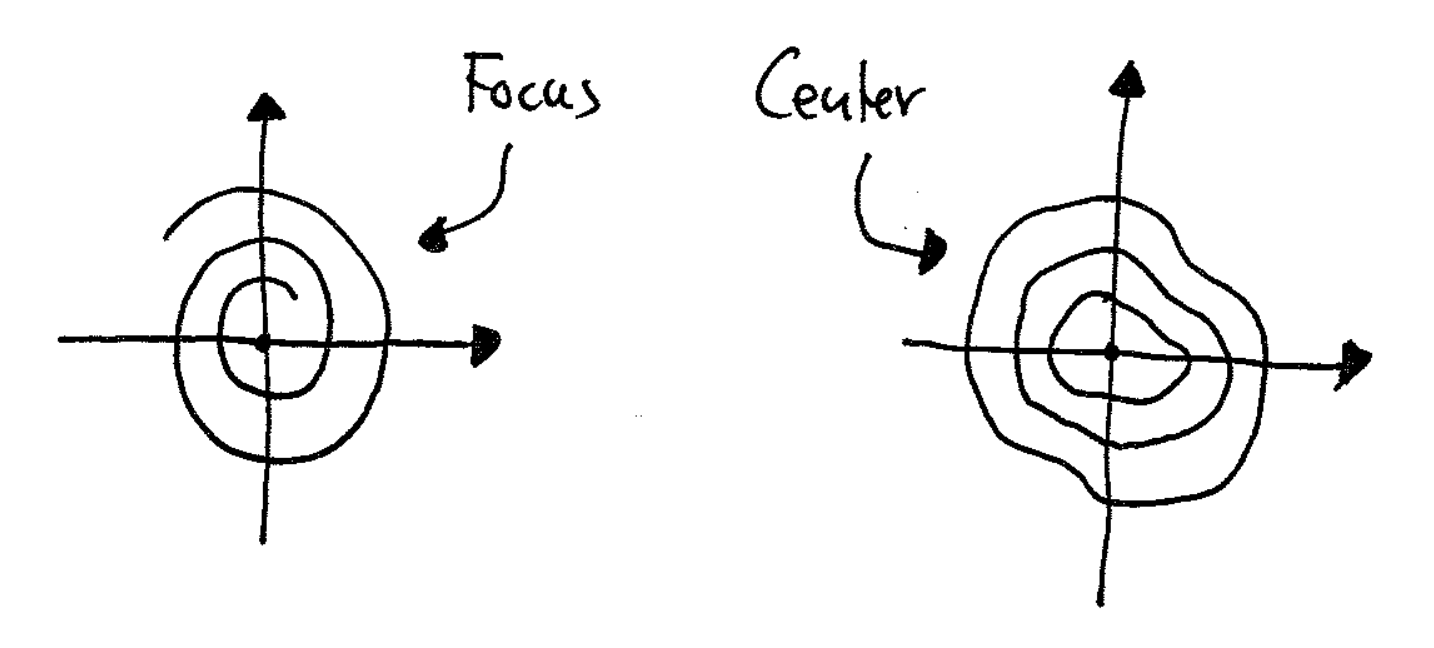}
\caption{A focus and a center.}
\label{fFocusCenter}
\end{figure}

\begin{rem} \label{rBlackBox}
The $\varepsilon$-method is most useful if there exists a fast algorithm for evaluating the polynomials of interest. The determinant of an $n \times n$ matrix for example has $n!$ terms, so the time to evaluate it directly is proportional to $n!$. If we use Gauss elimination on the matrix first, the time needed drops to $n^3$.
\end{rem}

For the remainder of this section we will look at an application of these methods to the Poincar\'e center problem. We start by considering the well known system of differential equations
\begin{align*}
	\dot{x} &= -y \\
	\dot{y} &= x
\end{align*}
whose integral curves are circles around the origin. Let's now disturb these equations with polynomials $P$ and $Q$ whose terms have degree at least $2$:
\begin{align*}
	\dot{x} &= -y + P\\
	\dot{y} &= x + Q.
\end{align*}
Near zero the integral curves of the disturbed system are either closed or not. In the second case one says that the equations have a \term{focus} in $(0,0)$ while in the first case they have a \term{center} (see Figure \ref{fFocusCenter}).

The condition of having a center is closed in the space of all $(P,Q)$:

\begin{thm}[Poincar\'e] \label{tDeg2}
There exists an infinite series of polynomials $f_i$ in the coefficients of $P$ and $Q$ such that
\[
\text{
$\begin{matrix}
	\dot{x} = -y + P\\
	\dot{y} = x + Q
\end{matrix}\,$
has a center $\iff f_i(P,Q) = 0$ for all $i$.}
\]
We call $f_i(P,Q)$ the $i$-th \term{focal value} of $(P,Q)$. 
%\end{thm}

%\begin{rem}
If the terms of $P$ and $Q$ have degree at most $d$ then the $f_i$ describe an algebraic variety $X_\infty$ in the finite-dimensional space of pairs $(P,Q)$. This variety is called the \term{center variety}. 
\end{thm}

\begin{rem}
By Hilbert's Basis Theorem $I_\infty := (f_0,f_1,\dots)$ is finitely generated. Unfortunately, Hilbert's Basis Theorem is not constructive, so it is a priory unknown how many generators $I_\infty$ has. It is therefore useful to consider the $i$-th partial center varieties $X_i = V(f_0,\dots,f_i)$.
\end{rem}

\begin{figure} 
\includegraphics*[width=12cm]{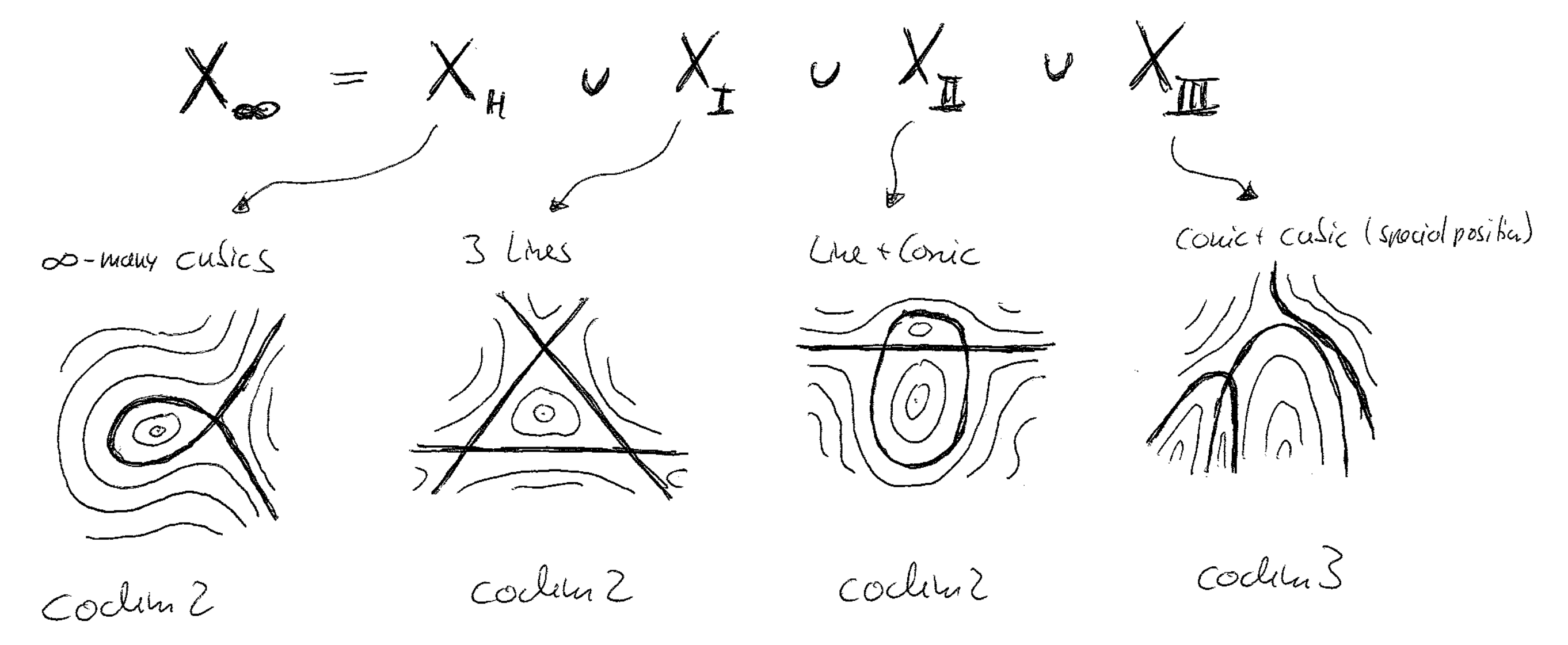}
\caption{Geometric interpretation of the components of the center variety in the case $d=2$.}
\label{fDeg2Cases}
\end{figure}

The following is known:

\begin{thm}
If $d=2$ then the center variety has four components
$$
X_\infty = X_H \cup X_{III} \cup X_{II} \cup X_{I} \subset \AZ^6,
$$
three of codimension $2$ and one of codimension $3$. Moreover $X_\infty = X_3$.
\end{thm}

\begin{proof}
Decompose $I_3=(f_1,f_2,f_3)$ with a computer algebra system and show that all solutions
do have a center \cite{frommer}, \cite{schlomiukTransactions}.
\end{proof}

Looking at algebraic integral curves one even obtains a geometric interpretation of the components in this case (see Figure \ref{fDeg2Cases}).

\newcommand{\Zoladec}{Zoladec}

For $d=3$ almost nothing is known. The best results so far are lists of centers given by \Zoladec \cite{ZoladekRational}, \cite{ZoladekCorrection}. The problem from a computer algebra perspective is that the $f_i$ are too large to be handled, already $f_5$ has $5348$ terms and it is known that $X_\infty \not= X_{10}$. 

\begin{figure}
\includegraphics*[width=12cm]{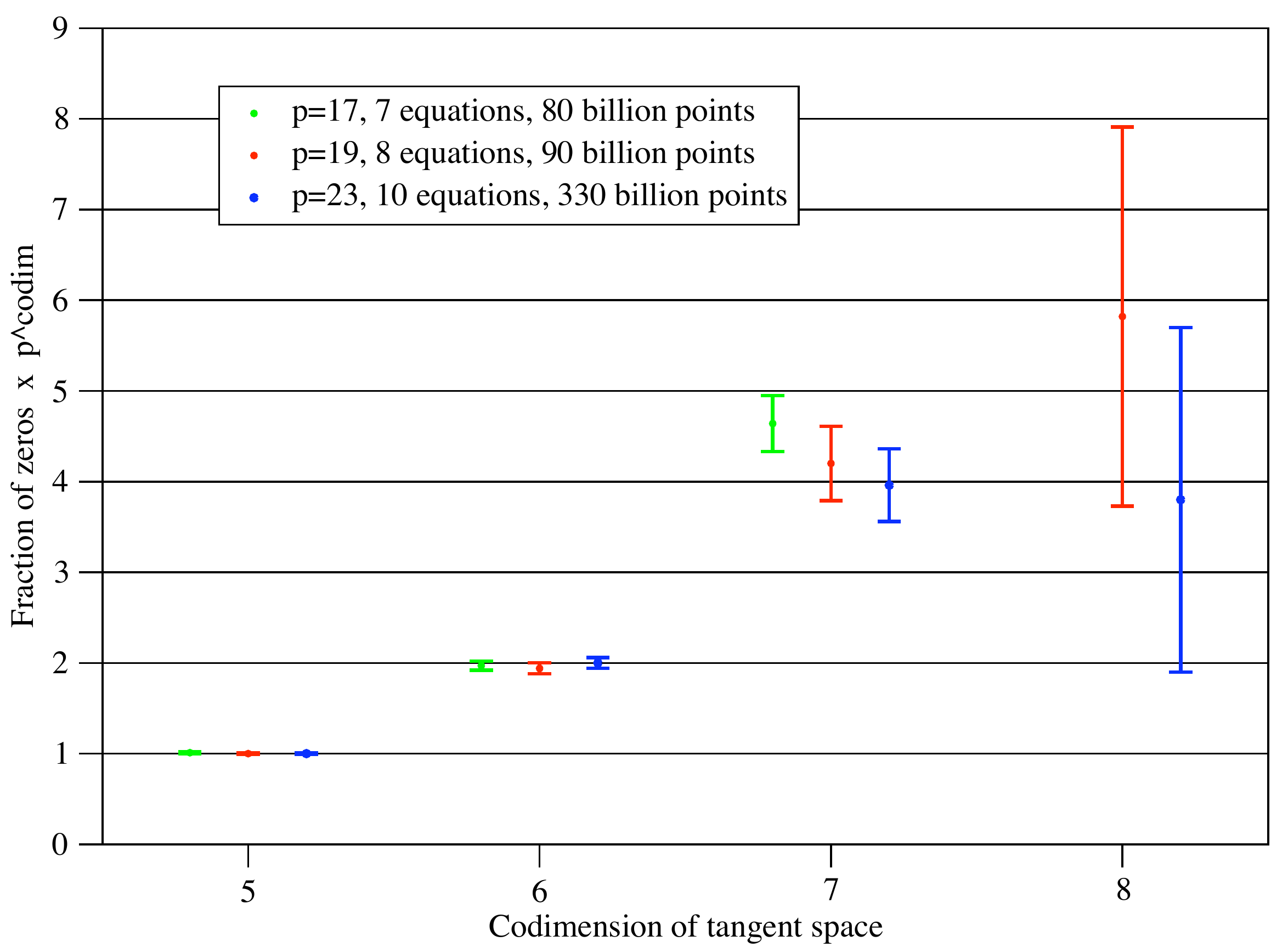}
\caption{Measurements for the Poincar\'e center problem with $d=3$.}
\label{fComponentsCut}
\end{figure}

\begin{experiment}
Fortunately for our method, Frommer \cite{frommer} has devised an algorithm to calculate $f_i(P,Q)$ for given $(P,Q)$. A closer inspection shows that Frommer's Algorithm works over finite fields and will also calculate $f_i(P+\varepsilon P',Q+\varepsilon Q')$. So we have all ingredients to use Heuristic \ref{hTangents}. Using a fast C++ implementation of Frommer's Algorithm by Martin Cremer and Jacob Kr\"oker \cite{strudelweb} we first check our method on the known degree $2$ case. For this we evaluate $f_1,\dots,f_{10}$  for $d=2$ at $1.000.000$ random points in characteristic $23$. This gives
\begin{verbatim}
  codim tangent space = 0: 5
  codim tangent space = 1: 162
  codim tangent space = 2: 5438
  codim tangent space = 3: 88
\end{verbatim}
Heuristic \ref{hTangents} translates this into 
\begin{verbatim}
  codim 0 components: 0.00 +/- 0.00
  codim 1 components: 0.00 +/- 0.00
  codim 2 components: 2.87 +/- 0.10
  codim 3 components: 1.07 +/- 0.29
\end{verbatim}
This agrees well with Theorem \ref{tDeg2}. 

For $d=3$ we obtain the measurements in Figure \ref{fComponentsCut}. 
One can check these results against \Zoladec's lists as depicted in Table \ref{tList}. Here the measurements agree in codimension $5$ and $6$. In codimension $7$ there seem to be 8 known families while we only measure $4$. Closer inspection of the known families reveals that $CR_5$ and $CR_7$ are contained in $CD_4$ and that $CR_{12}$ and $CR_{16}$ are contained in $CD_2$ \cite{zentrum}. After accounting for this our measurement agrees with \Zoladec's results and we conjecture that \Zoladec's lists are complete up to codimension $7$. 
\end{experiment}

\begin{table}
\caption{Known families of cubic centers that could have codimension below $8$ in $\AZ^{14}$ \cite{ZoladekCorrection}.}
\label{tList}

\begin{tabular}{|c|c|c|} \hline
Type & Name & Codimension \\ \hline  
Darboux & $CD_1$ & 5 \\
Darboux & $CD_2$& 6 \\
Darboux & $CD_3$ & 7 \\
Darboux & $CD_4$ & 7 \\
Darboux & $CD_5$ & 7 \\ 
\hline
Reversible & $CR_1$   & $\ge 6$ \\
Reversible & $CR_5$   & $\ge 7$ \\
Reversible & $CR_7$   & $\ge 7$ \\
Reversible & $CR_{11}$ & $\ge 7$ \\
Reversible & $CR_{12}$ & $\ge 7$ \\ 
Reversible & $CR_{16}$ & $\ge 7$ \\ 
\hline
\end{tabular}

\end{table}

%%%%%%%%%%%%%%%%%%%%%%%%%%%%%%%%%%%%%%%%%%%%%
\section{Existence of a Lift to Characteristic Zero} \label{sExLift}
%%%%%%%%%%%%%%%%%%%%%%%%%%%%%%%%%%%%%%%%%%%%%
\nosubsections

Often one is not interested in characteristic $p$ solutions, but in solutions over $\CC$. Unfortunately,
not all solutions over $\FF_p$ lift to characteristic $0$. 

\begin{figure}[h!]
\includegraphics*[width=10cm]{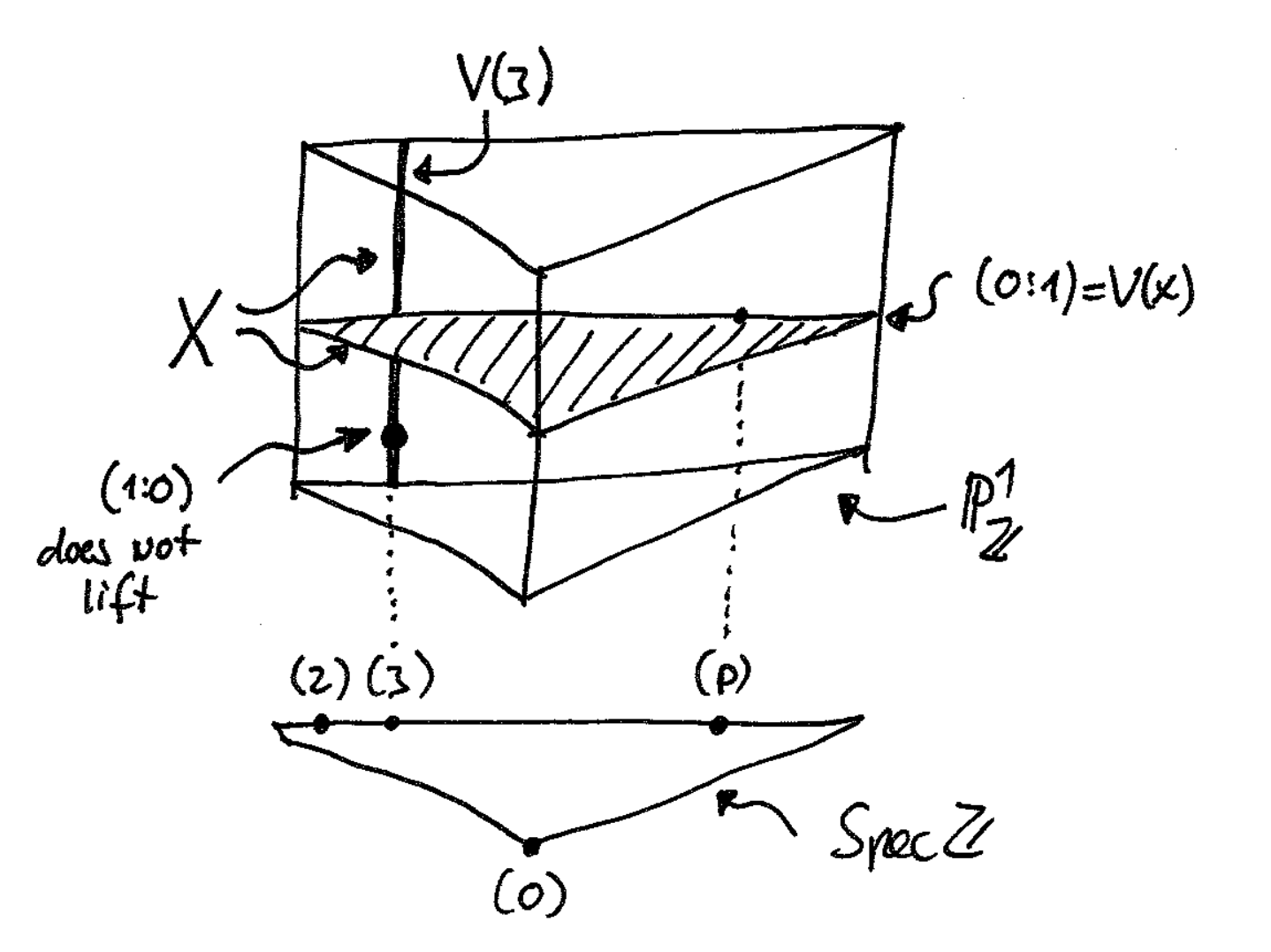}
\caption{The vanishing set of $3x$ in $\PP^1_\ZZ$ over $\spec \ZZ$}
\label{fV3x}
\end{figure}

\begin{example}
Consider the variety $X = V(3x) \subset \PP^1_\ZZ$ over $\spec \ZZ$. As depicted in Figure \ref{fV3x}, $X$ decomposes into two components: $V(3) = \PP^1_{\FF_3}$ which lives only over $\FF_3$ and $V(x)= \{(0:1)\}$ which has fibers over all of $\spec \ZZ$. In particular, the point $(1:0) \in \PP^1_{\FF_3} \subset X$ does not lift to characteristic $0$.
\end{example}

To prove that a given solution point over $\FF_p$ does lift to characteristic zero the following tool is very helpful:

\newcommand{\Xp}{X_{\FF_p}}
\newcommand{\Yp}{Y_{\FF_p}}
\newcommand{\Zp}{Z_{\FF_p}}
\newcommand{\XZZ}{X_{\ZZ}}
\newcommand{\YZZ}{Y_{\ZZ}}
\newcommand{\ZZZ}{Z_{\ZZ}}

\begin{prop}[Existence of a Lifting] \label{pExLift}
Let $X \subset Y \subset \AZ^n_\ZZ$ be varieties with $\dim \Yp = \dim \YZZ -1$ for all $p$ and $X \subset Y$ determinantal, i.e.
there exists a vector bundle morphism
\[
	\phi \colon E \to F
\]
on $Y$ and a number $r \le \min(\rank E, \rank F)$ such that $X = X_r(\phi)$ is the locus where $\phi$ has rank at most $r$. If $x \in \Xp$ is a point with
\[
	\dim T_{\Xp,x} = \dim \Yp - (\rank E-r)(\rank F -r)
\]
then $X$ is smooth in $x$ and there exists a component $Z$ of $\XZZ$ containing $x$ and having
a nonzero fiber over $(0)$.
\end{prop}

\begin{proof}
Set $d = \dim \Yp - (\rank E-r)(\rank F -r)$.
Since $\Xp$ is determinantal, we have
$$\dim \Zp \ge d$$
for every irreducible component $\Zp$ of $\Xp$ and $d$ is the expected dimension of $\Zp$  \cite[Ex.\,10.9, p.\,245]{Ei95}. If $\Zp$ contains the point $x$ we obtain
\[
	d \le \dim \Zp \le \dim T_{\Zp,x} \le \dim T_{\Xp,x} = d
\]
by our assumptions. So $\Zp$ is of dimension $d$ and smooth in $x$. Let now $\ZZZ$ be a component
of $\XZZ$ that contains $\Zp$ and $x$. Since $\XZZ$ is determinantal in $\YZZ$ and $\dim \YZZ = \dim \Yp + 1$ we have
\[
	\dim \ZZZ \ge d+1.
\]
Since $\dim \Zp = d$ the fiber of $\ZZZ$ over $p$ cannot contain all of $\ZZZ$. Indeed, in this case we would have $\Zp = \ZZZ$ since both are irreducible, but $\dim \Zp \not= \dim \ZZZ$. It follows
that $\ZZZ$ has nonempty fibers over an open subset of $\spec \ZZZ$ and therefore also over
$(0)$ \cite{smallFields}, \cite{newfamily}.
\end{proof}

\begin{figure}
\includegraphics*[width=11cm]{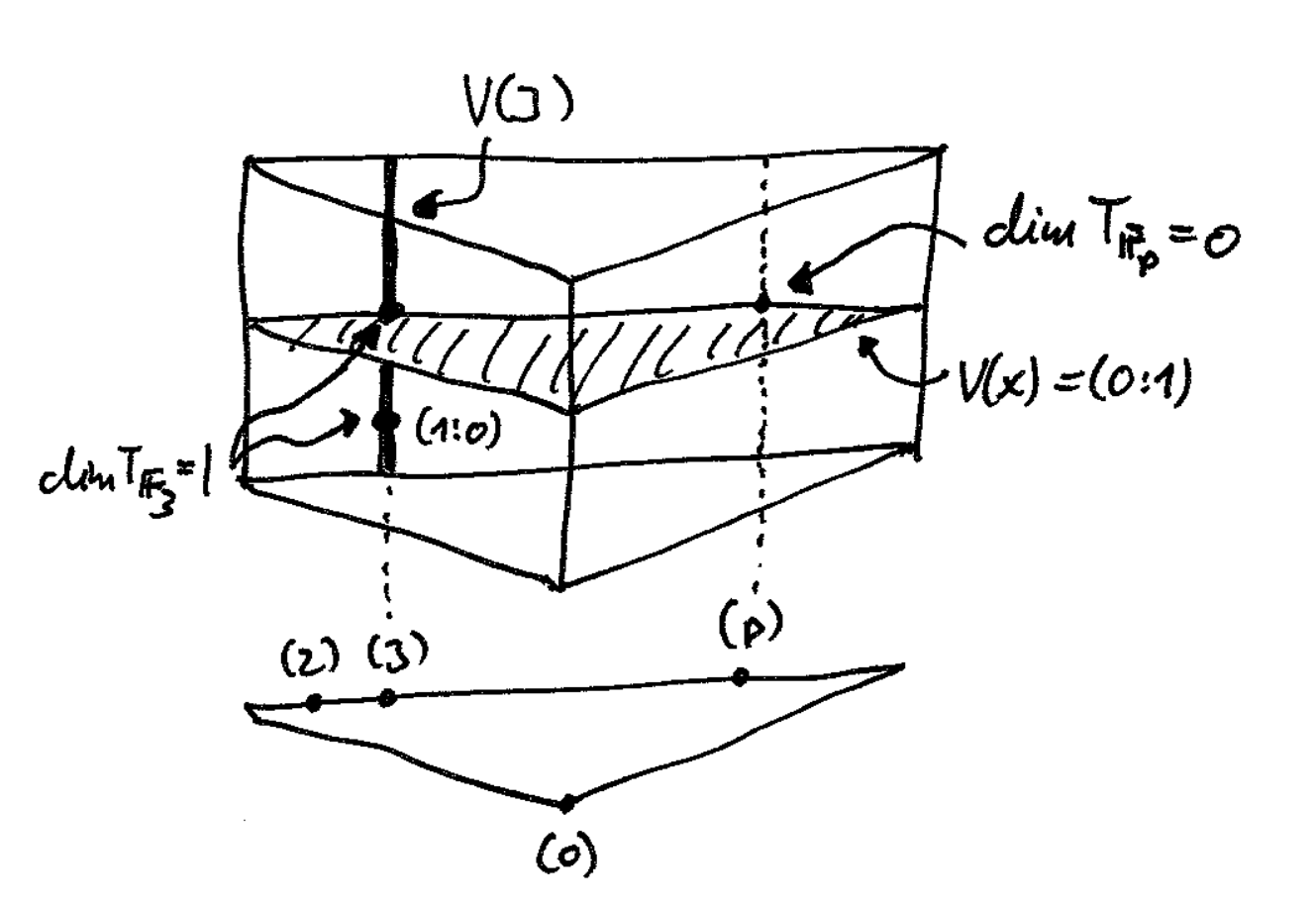}
\caption{Tangent spaces in several points of $X=V(3x) \subset \PP^1_\ZZ$ over $\spec \ZZ$}
\label{fV3xTangent}
\end{figure}

\begin{example}
The variety $X=V(3x)$ is determinantal on $Y=\PP^1_\ZZ$ since it is the rank $0$ locus of the 
vector bundle morphism
\[	
	 \phi \colon \sO_{\PP^1_\ZZ} \xrightarrow{3x} \sO_{\PP^1_{\ZZ}}(1).
\]
Furthermore $\dim \PP^1_{\FF_p} = 1 = \dim \PP^1_\ZZ -1$ for all $p$. The expected dimension of $\Xp$ is therefore $1 - (1-0)\cdot (1-0) = 0$. As depicted in Figure \ref{fV3xTangent}
we have three typical examples:
\begin{enumerate}
\item $x = (0:1)$ over $\FF_p$ with $p \not=3$. Here the tangent space over $\FF_p$ is zero dimensional and the point lifts according to Proposition \ref{pExLift}.
\item $x = (0:1)$ over $\FF_3$. Here the tangent space is $1$-dimensional and Proposition \ref{pExLift} does not apply. Even though the point does lift.
\item $x = (1:0)$ over $\FF_3$. Here the tangent space is also $1$-dimensional and Proposition
\ref{pExLift} does not apply. In this case the point does not lift.
\end{enumerate}
\end{example}

%\begin{rem}
%This is a generalization of Hensel-lifting.
%\end{rem}

This method has been used first by Frank Schreyer  \cite{smallFields} to construct new surfaces in $\PP^4$ which are not of general type. The study of such surfaces started in $1989$ when Ellingsrud and Peskine showed that their degree is bounded  
\cite{EllingsrudPeskine} and therefore only finitely many families exist. Since then the degree bound has been sharpened by various
authors, most recently by \cite{DeckerSchreyer} to $52$. On the other hand a classification is only known up to degree $10$ and examples are known up to degree $15$ (see \cite{DeckerSchreyer} for an overview and references). 

Here I will explain how Cord Erdenberger, Katharina Ludwig and I found a new family of rational surfaces $S$ of degree $11$ and sectional genus $11$ in $\PP^4$ with finite field experiments.

Our plan is to realize $S$ as a blowup of $\PP^2$. First we consider some restrictions on the linear system that embeds $S$ into $\PP^4$:

\begin{prop} \label{p-formulae1}
Let $S = \PZ^2_\CC(p_1,\dots,p_l)$ be the blowup of $\PZ^2_\CC$ in $l$ distinct points. We denote by $E_1,\dots,E_l$ the corresponding exceptional divisors and by $L$ the pullback of a general line in $\PP^2_\CC$ to $S$. Let
$
    | aL - \sum_{i=1}^l b_i E_i |
$
be a very ample linear system of dimension four and set $\beta_j = \#\{i \suchthat b_i=j\}$. Then
\begin{align*}
    d     &= a^2 - \sum_j \beta_j j^2\\
    \pi   &=  { a-1 \choose 2} - \sum_{j}\beta_j { j\choose 2}\\
    K^2 &= 9 - \sum_j \beta_j.
\end{align*}
where $d$ is the degree, $\pi$ the sectional genus and $K$ the canonical divisor of $S$.
\end{prop}

\begin{proof} Intersection theory on $S$ \cite[Corollary 4.1]{newfamily}. \end{proof}

By the double point formula for surfaces in $\PP^4$  \cite[Appendix A, Example 4.1.3]{Ha} a rational surface of degree $11$ and sectional genus $11$ must satisfy $K^2 = -11$. For fixed $a$ the equations above can be solved by integer programming, using for example the algorithm described in Chapter 8 of \cite{CLOusingAG}. 

In the case $a <9$ we find that there are no solutions. For $a=9$ the only solution is $\beta_3 = 1$, $\beta_2 = 14$ and $\beta_1 = 5$. Our first goal is therefore to find $5$ simple points, $14$ double points and one triple point in $\PP^2$ such that the ideal of the union of these points contains $5$ polynomials of degree $9$.

To make the search fast, we would like to use characteristic $2$. The difficulty here is that $\PP^2$ contains only $7$ rational points, while we need $20$. Our solution to this problem was to choose
\begin{align*}
	 P &\in \PP^2(\FF_2) &
	 Q &\in \PP^2(\FF_{2^{14}}) &
	 R &\in \PP^2(\FF_{2^5})
\end{align*}
such that the Frobenius orbit of $Q$ and $R$ are of length $14$ and $5$ respectively. The ideals of the orbits are then defined over $\FF_2$. 
\begin{verbatim}
  -- define coordinate ring of P^2 over F_2
  F2 = GF(2)
  S2  = F2[x,y,z]

  -- define coordinate ring of P^2 over F_2^14 and F_2^5
  St  = F2[x,y,z,t]
  use St; I14 = ideal(t^14+t^13+t^11+t^10+t^8+t^6+t^4+t+1); S14 = St/I14
  use St; I5 = ideal(t^5+t^3+t^2+t+1); S5 = St/I5

  -- the random points
  use S2; P = matrix{{0_S2, 0_S2, 1_S2}}
  use S14;Q = matrix{{t^(random(2^14-1)), t^(random(2^14-1)), 1_S14}} 
  use S5; R = matrix{{t^(random 31), t^(random 31), 1_S5}}

  -- their ideals
  IP = ideal ((vars S2)*syz P)
  IQ = ideal ((vars S14)_{0..2}*syz Q)
  IR = ideal ((vars S5)_{0..2}*syz R)

  -- their orbits
  f14 = map(S14/IQ,S2); Qorbit = ker f14
  degree Qorbit   -- hopefully degree = 14

  f5 = map(S5/IR,S2); Rorbit = ker f5
  degree Rorbit   -- hopefully degree = 5
\end{verbatim}
If $Q$ and $R$ have the correct orbit length we calculate $|9H-3P-2Q-R|$
\begin{verbatim}
  -- ideal of 3P
  P3 = IP^3;

  -- orbit of 2Q
  f14square = map(S14/IQ^2,S2); Q2orbit = ker f14square;

  -- ideal of 3P + 2Qorbit + 1Rorbit
  I = intersect(P3,Q2orbit,Rorbit);

  -- extract 9-tics
  H = super basis(9,I)
  rank source H   -- hopefully affine dimension = 5
\end{verbatim}
If at this point we find $5$ sections, we check that there are no unassigned base points
\begin{verbatim}
  -- count basepoints (with multiplicities)
  degree ideal H   -- hopefully degree = 1x6+14x3+1x5 = 53
\end{verbatim}
If this is the case, the next difficulty is to check if the corresponding linear system is very ample. On the one hand this is an open condition, so it should be satisfied by most examples, on the other hand we are in characteristic $2$, so exceptional loci can have very many points. An irreducible divisor for example already contains approximately half of the rational points. 
\begin{verbatim}
  -- construct map to P^4
  T = F2[x0,x1,x2,x3,x4]
  fH = map(S2,T,H);

  -- calculate the ideal of the image
  Isurface = ker fH; 

  -- check invariants
  betti res coker gens Isurface
  codim Isurface    -- codim = 2
  degree Isurface   -- degree = 11
  genera Isurface   -- genera = {0,11,10}

  -- check smoothness
  J = jacobian Isurface;
  mJ = minors(2,J) + Isurface;
  codim mJ  -- hopefully codim = 5
\end{verbatim}
Indeed, after about $100.000$ trials one comes up with the points
\begin{verbatim}
  use S14;Q = matrix{{t^11898, t^137, 1_S14}}
  use S5; R = matrix{{t^6, t^15, 1_S5}}
\end{verbatim}
\newcommand{\PtwoZZ}{{\PZ^2_\ZZ}}
\newcommand{\sOPtwoZZ}{{\sO_\PtwoZZ}}
These satisfy all of the above conditions and prove that rational surfaces of degree $11$ and sectional genus $11$ in $\PP^4$ exist in over $\FF_2$. 

As a last step we have to show that this example lifts to char $0$. For this we consider
the morphism
\[
    \tau_k \colon H^0(\sOPtwoZZ(a)) \to \sOPtwoZZ(a) \oplus
                                3 \sOPtwoZZ(a-1) \oplus
                                \dots \oplus
                                { k+2 \choose 2} \sOPtwoZZ(a-k)
\]
on $\PtwoZZ$ that associates to each polynomial of degree $a$ the coefficients of its Taylor expansion up to degree $k$ in a given point $P$.

\begin{lem}
If $a > k$ then the image of $\tau_k$ is a vector bundle $\sF_k$ of rank ${k+2 \choose 2}$ over $\spec \ZZ$.
\end{lem}

\begin{proof}
In each point we consider an affine $2$-dimensional neighborhood where we can choose the
${k+2 \choose 2}$ coefficients of the affine Taylor expansion independently. This shows that
the image has at least this rank everywhere. If follows from the Euler relation for homogeneous polynomials
$$
	x\frac{df}{dx}+y\frac{df}{dy}+z\frac{df}{dz} = (\deg f)\cdot f
$$
 that this is also the maximal rank. 
\end{proof}

Now set $Y_\ZZ = \Hilb_{1,\ZZ} \times \Hilb_{14,\ZZ} \times \Hilb_{5,\ZZ}$ where $\Hilb_{k,\ZZ}$ denotes the Hilbert scheme of $k$ points in $\PP^2_\ZZ$ over $\spec \ZZ$, and let
\[
    X_\ZZ = \{ (p,q,r) \suchthat h^0(9L - 3p - 2q - 1r) \ge 5\} \subset Y_\ZZ
\]
be the subset where the linear system of nine-tics with a triple point in $p$, double points in $q$ and
single base points in $r$ is at least of projective dimension $4$.

\begin{prop}
There exist vector bundles $E$ and $F$ of ranks $55$ and $53$ respectively on $Y_\ZZ$ and a morphism
\[
    \phi  \colon E \to F
\]
such that $X_{50}(\phi) = X_\ZZ$.
\end{prop}

\begin{proof}
On the Cartesian product
\xycenter{
     \Hilb_{d,\ZZ} \times \PtwoZZ \ar[r]^-{\pi_2} \ar[d]^-{\pi_1} & \PtwoZZ \\
      \Hilb_{d,\ZZ}
            }
we have the morphisms
\[
    \pi_2^*\tau_k \colon H^0(\sOPtwoZZ(9))\otimes \sO_{\Hilb_{d,\ZZ} \times \PtwoZZ} \to \pi_2^* \sF_k.
\]
Let now $P_d \subset  \Hilb_{d,\ZZ} \times \PtwoZZ$ be the universal set of points. Then
$P_d$ is a flat family of degree $d$ over $\Hilb_{d,\ZZ}$ and
\[
    \sG_{k} := (\pi_1)_* ((\pi_2^*\sF_k)|_{P_d})
\]
is a vector bundle of rank $d{ k+2 \choose 2}$ over $\Hilb_{d,\ZZ}$. On
\[
      Y_\ZZ = \Hilb_{1,\ZZ} \times \Hilb_{14,\ZZ} \times \Hilb_{5,\ZZ}
\]
the induced map
\[
    \phi \colon H^0(\sOPtwoZZ(9))\otimes \sO_{X_\ZZ} \xrightarrow{\tau_2\oplus\tau_1\oplus\tau_0}
                         \sigma_1^* \sG_{2} \oplus \sigma_{14} ^* \sG_{1}
                         \oplus \sigma_{5}^* \sG_{0}
\]
has the desired properties, where $\sigma_d$ denotes the projection to $\Hilb_{d,\ZZ}$.
\end{proof}

So we have to show that the tangent space of $X_{\FF_2}$ in our base locus has codimension
$(55-50)(53-50)=15$. This can be done by explicitly calculating the differential of $\phi$ in our given
base scheme using the $\varepsilon$-method. The script is too long for this paper, but can be downloaded at \cite{ratsurfweb}. Indeed, we find that the codimension of the tangent space is $15$, so this shows that our example lies on an irreducible component that is defined over an open subset of $\spec \ZZ$. 

\begin{rem}
The overall time to find smooth surfaces that lift to characteristic zero can be substantially reduced if one calculates the tangent space of a given point $(P,Q,R)$ in the Hilbert scheme $X_\ZZ$ directly after establishing $|9H-3P-2Q-R| = \PP^4$. One then needs to check very ampleness only for smooth points of $X_\ZZ$. This is useful since the tangent space calculation is just a linear question, while the check for very ampleness requires Gr\"obner bases. We use a very fast $C$-implementation by Jakob Kr\"oker to do the whole search algorithm up to checking smoothness. Only the (very few) remaining examples are then checked for very ampleness using Macaulay $2$.
\end{rem}

\begin{figure}
\includegraphics*[width=8cm]{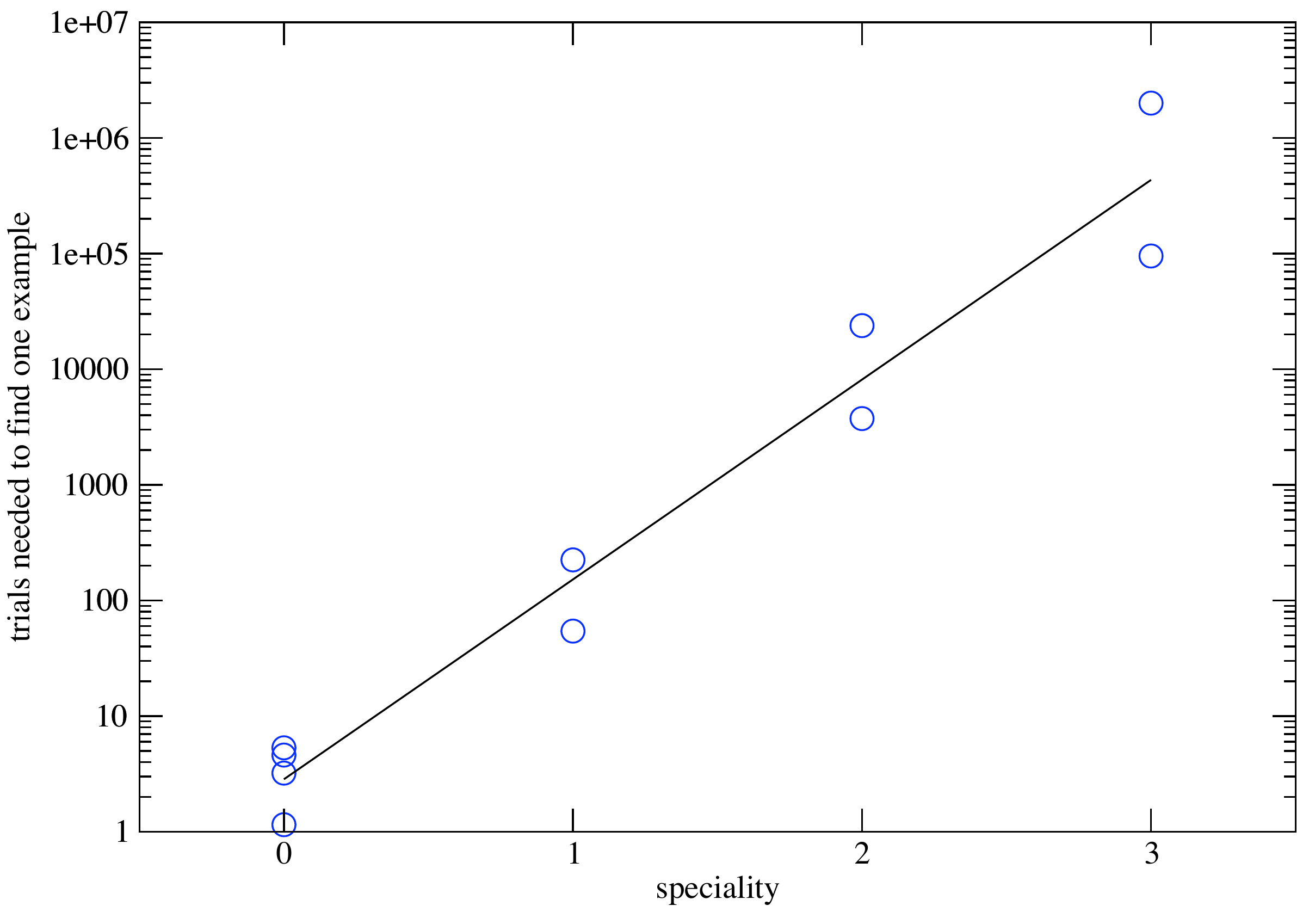}
\caption{The difficulty of finding a surface grows exponentially with the speciality}
\label{fTrials}
\end{figure}

\begin{experiment}
We also tried to reconstruct the other known rational surfaces in $\PP^4$ with our program. The number
of trials needed is depicted in Figure \ref{fTrials}. The expected codimension of $X$ in the corresponding
Hilbert scheme turns out to be $5$ times the speciality $h^1(\sO_X(1))$ of the surface. As expected, the logarithm of the number of trials needed to find a surface is proportional to the codimension of $X$.
\end{experiment}

\begin{rem}
We could not reconstruct all known families. The reason for this is that we only look at examples where the base points of a given multiplicity form an irreducible Frobenius orbit. In some cases such examples do not exist for geometric reasons.
\end{rem}

\begin{experiment}
Looking at the linear system $|14H - 4P - 3Q-R|$ with $\deg P = 8$, $\deg Q = 6$ and $\deg R =2$, we 
find rational surfaces of degree $12$ and sectional genus $12$ in $\PP^4$ with this method (not published) . 
\end{experiment}

%%%%%%%%%%%%%%%%%%%%%%%%%%%%%%%%%%%%%%%%%%%%%
\section{Finding a Lift} \label{sFindLift}
%%%%%%%%%%%%%%%%%%%%%%%%%%%%%%%%%%%%%%%%%%%%%
\nosubsections
%\suppressfloats 
In some good cases characteristic $p$ methods even allow one to find a solution over $\QQ$ quickly.
Basically this happens when the solution set is zero dimensional with two different flavors. 

\begin{figure}[h]
\includegraphics*[width=10cm]{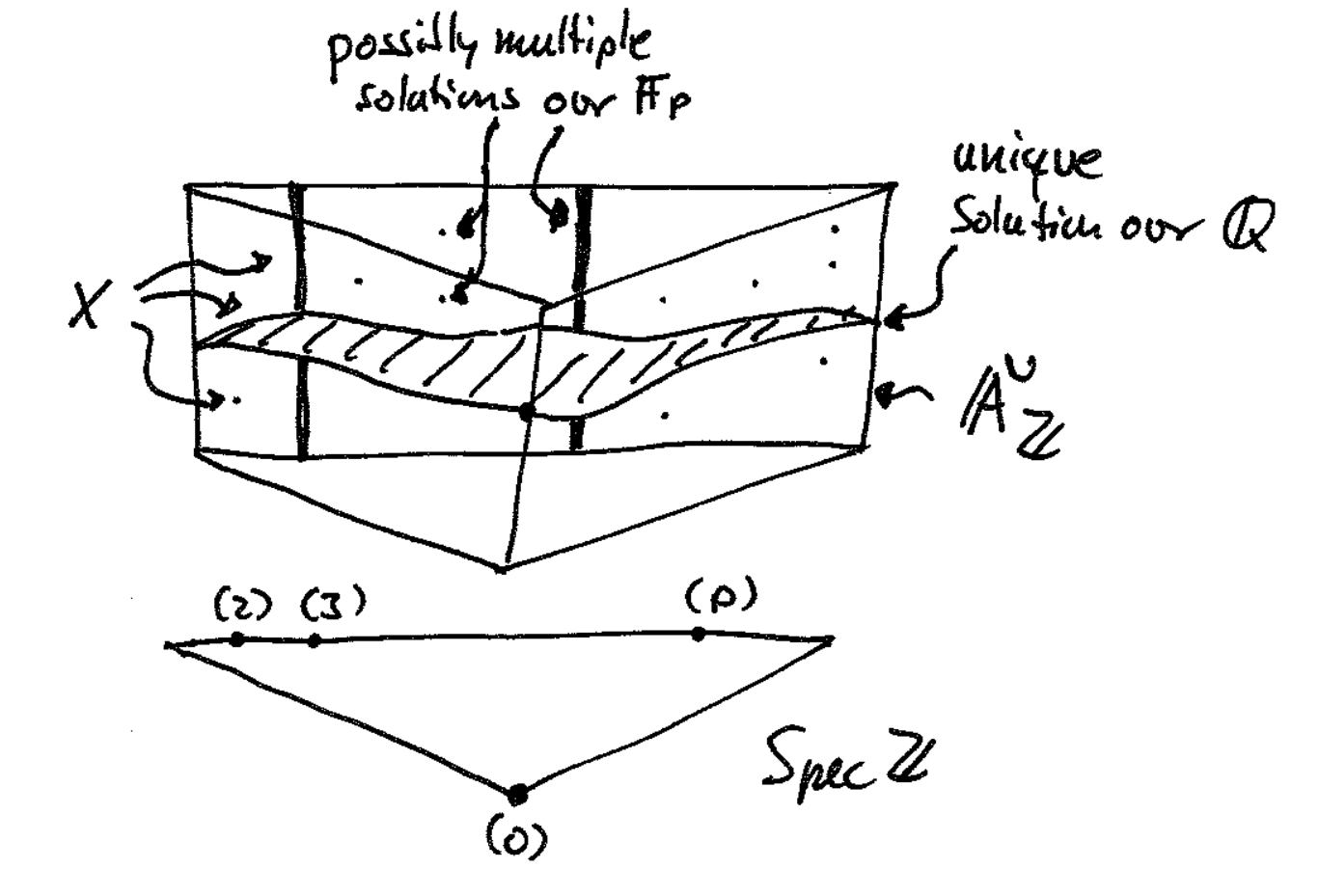}
\caption{A scheme over $\spec \ZZ$ with a {\sl unique} solution over $\bar{\QQ}$, possibly with high multiplicity}
\label{fUniqueQQ}
\end{figure}

The first good situation, depicted in Figure \ref{fUniqueQQ},  arises when $X = V(f_1,\dots, f_m) \subset \AZ^n_{\ZZ}$ has a unique solution over $\bar{\QQ}$, maybe with high multiplicity. In this case it follows that the solution is defined over $\QQ$.

\begin{alg} \label{aLiftZZ}
If the coordinates of the unique solution over $\QQ$ are even in $\ZZ$ one can find this solution 
as follows:
\begin{enumerate}
\item Reduce mod $p_i$ and test {\sl all} points in $\FF_{p_i}^n$
\item Find many primes $p_i$ with a unique solution in $\FF_{p_i}^n$
\item Use Chinese remaindering to find a solution mod $\prod_i p_i >> 0$.
\item Test if this is a solution over $\ZZ$. If not, find more primes $p_i$ with
unique solutions over $\FF_{p_i}$. 
\end{enumerate}
\end{alg}
 
\begin{rem}
Even if the solution $y$ over $\QQ$ is unique, there can be several solutions over $\FF_{p}$. Since the
codimension of points in $\AZ^n$ is $n$ we expect that the probability of a random point $x \in \AZ^n$
to satisfy $x \in X$ is $\frac{1}{p^n}$ by Heuristic \ref{hHigherCodim}. We therefore expect that the probability of $x \not\in X$
for all points $x \not=y$ is
\[
	\left(1-\frac{1}{p^n}\right)^{p^n-1} \approx \frac{1}{e}.
\]
\end{rem}

\begin{experiment} \label{eLiftZZ}
Let's use this Algorithm \ref{aLiftZZ} to solve
\begin{align*}
	-8 x^{2}-x y-7 y^{2}+5238 x-11582y-7696 &= 0 \\
	4 x y-10 y^{2}-2313 x-16372 y-6462 &= 0
\end{align*}
For this
we need a function that looks at all points over a given prime:
\begin{verbatim}
  allPoints = (I,p) -> (
     K = ZZ/p;
     flatten apply(p,i->
     	  flatten apply(p,j->
	           if (0==codim sub(I,matrix{{i*1_K,j*1_K}})) 
	           then {(i,j)}
	           else {}
     	  ))
     )
\end{verbatim}
With this we look for solutions of our equations over the first nine primes.
\begin{verbatim}
  R = ZZ[x,y]
  -- the equations
  I = ideal (-8*x^2-x*y-7*y^2+5238*x-11582*y-7696,
     4*x*y-10*y^2-2313*x-16372*y-6462)      
      
  -- look for solutions
  tally apply({2,3,5,7,11,13,17,19,23},p->(p,time allPoints(I,p)))
\end{verbatim}
We obtain:
\begin{verbatim}
  o8 = Tally{(2, {(0, 0)}) => 1                             
           (3, {(0, 2), (1, 0), (2, 0)}) => 1
           (5, {(4, 1)}) => 1
           (7, {(2, 3), (5, 5)}) => 1
           (11, {(2, 7), (8, 1)}) => 1
           (13, {(3, 4), (12, 6)}) => 1
           (17, {(10, 8)}) => 1
           (19, {(1, 3), (1, 17), (18, 5), (18, 18)}) => 1
           (23, {(15, 8)}) => 1
\end{verbatim}
As expected for the intersection of two quadrics we find at most $4$ solutions. Over four primes we find unique solutions, which is reasonably close to the expected number $9/e \approx 3.31$. We now combine the information over these four primes using the Chinese remainder Theorem.
\begin{verbatim}
  -- Chinese remaindering
  -- given solutions mod m and n find
  -- a solution mod m*n
  -- sol1 = (n,solution)
  -- sol2 = (m,solution)
  chinesePair = (sol1, sol2) -> (
       n = sol1#0;an = sol1#1;
       m = sol2#0;am = sol2#1;
       drs = gcdCoefficients(n,m);
       -- returns {d,r,s} so that a*r + b*s is the
       -- greatest common divisor d of a and b.
       r = drs#1;
       s = drs#2;
       amn = s*m*an+r*n*am;
       amn = amn - (round(amn/(m*n)))*(m*n);
       if (drs#0) == 1 then (m*n,amn) else print "m and n not coprime"
       )

  -- take a list {(n_1,s_1),...,(n_k,s_k)}
  -- and return (n,a) such that
  -- n = n_1 * ... * n_k   and
  -- s_i = a mod n_i
  chineseList = (L) -> (fold(L,chinesePair))

  -- x coordinate
  chineseList({(2,0),(5,4),(17,10),(23,15)})
  -- y coordinate
  chineseList({(2,0),(5,1),(17,8),(23,8)})
\end{verbatim}
This gives
\begin{verbatim}
  o11 = (3910, 1234)
  o12 = (3910, -774)
\end{verbatim}
i.e $(1234,-774)$ is the unique solution mod $3910 = 2\cdot5\cdot17\cdot23$. Substituting
this into the original equations over $\ZZ$ shows that this is indeed a solution over $\ZZ$.
\begin{verbatim}
  sub(I,matrix{{1234,-774}})

  o13 = ideal (0, 0)
\end{verbatim}
\end{experiment}

If the unique solution does not have $\ZZ$ but $\QQ$ coordinates then one can find
the solution using the extended Euclidean Algorithm \cite[Section 5.10]{vzGathen}.

\begin{example}
Let's try to find a small solution to the equation
\[
	\frac{r}{s} \equiv 7\mod 37.
\]
Each solution satisfies 
\[
	r = 7s + 37 t
\]
with $s$ and $t$ in $\ZZ$. Using the extended Euclidean Algorithm 
\begin{center}
\begin{tabular}{|c|c|c|c|c|}
\hline
& $r$ & $s$ & $t$ & $r/s$ \\
\hline
     & $37$ & $0$ & $1$ & \\
$-5$ &  $7$  & $1$ & $0$ & $7/1$ \\
$-3$ & $2$  & $-5$ & $1$ & $-2/5$ \\
     & $1$ & $16$ & $-3$ & $1/16$ \\
\hline
\end{tabular} 
\end{center}
we find the solution
\[
	1 = \gcd(7,37) = 7\cdot 16 + 37 \cdot(-3)
\]
to our linear equation. Observe, however, that the intermediate step in the
Euclidean Algorithm also gives solutions, most of them with small coefficients.
Indeed, $r/s=-2/5$ is a solution with $r,s \le \sqrt{37}$ which is the best that we can expect.
\end{example}

If we find a small solution by this method, we even can be sure that it is the only one satisfying the 
congruence:

\begin{prop}
There exist at most two solutions $(r,s)$ of 
$$r \equiv as+bt \mod m$$ 
that satisfy $r,s \le \sqrt{m}$.
If a solution satisfies $r,s \le \frac{1}{2}\sqrt{m}$, then this solution is unique.
\end{prop}

\begin{proof}
\cite[Section 5.10]{vzGathen}
\end{proof}

\begin{experiment}
Let's find a solution to 
\begin{align*}
176 x^{2}+148 x y+301 y^{2}-742 x+896y+768  &= 0 \\
-25 x y+430 y^{2}+33 x+1373 y+645 &= 0
\end{align*}
As in Experiment \ref{eLiftZZ} we search for primes with unique solutions
\begin{verbatim}
  I = ideal (176*x^2+148*x*y+301*y^2-742*x+896*y+768,
         -25*x*y+430*y^2+33*x+1373*y+645)
  tally apply({2,3,5,7,11,13,17,19,23,29,31,37,41},
         p->(p,time allPoints(I,p)))
\end{verbatim}
and obtain
\begin{verbatim}
  o10 = Tally{(2, {(1, 0)}) => 1                              
            (3, {(0, 0), (0, 1), (2, 0)}) => 1
            (5, {(3, 2), (4, 1)}) => 1
            (7, {(2, 6), (4, 0)}) => 1
            (11, {}) => 1
            (13, {(5, 10)}) => 1
            (17, {(5, 4), (9, 13), (11, 16), (12, 12)}) => 1
            (19, {(3, 15), (8, 6), (13, 15), (17, 1)}) => 1
            (23, {(15, 18), (19, 12)}) => 1
            (29, {(26, 15), (28, 9)}) => 1
            (31, {(7, 22)}) => 1
            (37, {(14, 18)}) => 1
            (41, {(0, 23)}) => 1
\end{verbatim}
Notice that there is no solution mod $11$. If there is a solution over $\QQ$ this means that
$11$ has to divide at least one of the denominators. Chinese remaindering gives a solution mod 
$2\cdot13\cdot31\cdot37\cdot41 = 1222702$:
\begin{verbatim}
  -- x coordinate
  chineseList({(2,1),(13,5),(31,7),(37,14),(41,0)})
  o11 = (1222702, 138949)
  -- y coordinate
  chineseList({(2,0),(13,10),(31,22),(37,18),(41,23)})
  o12 = (1222702, -526048)
\end{verbatim}
Substituting this into the original equations gives
\begin{verbatim}
  sub(I,matrix{{138949,-526048}})
  o13 = ideal (75874213835186, 120819022681578)
\end{verbatim}
so this is not a solution over $\ZZ$. To find a small possible solution over $\QQ$ we use an implementation of the extended Euclidean Algorithm from  \cite[Section 5.10]{vzGathen}.
\begin{verbatim}
  -- take (a,n) and calculate a solution to
  --   r =  as mod n
  -- such that r,s < sqrt(n).
  -- return (r/s)
  recoverQQ = (a,n) -> (
     r0:=a;s0:=1;t0:=0;
     r1:=n;s1:=0;t1:=1;
     r2:=0;s2:=0;t2:=0;
     k := round sqrt(r1*1.0);
     while k <= r1 do (
     	  q = r0//r1;
     	  r2 = r0-q*r1;
     	  s2 = s0-q*s1;
     	  t2 = t0-q*t1;
     	  --print(q,r2,s2,t2);
     	  r0=r1;s0=s1;t0=t1;
     	  r1=r2;s1=s2;t1=t2;
     	  );
     (r2/s2)
     )
\end{verbatim}
This yields
\begin{verbatim}
  -- x coordinate
  recoverQQ(138949,2*13*31*37*41)
  
        123
  o21 = ---
         22
\end{verbatim}
Notice that Macaulay reduced $246/44$ to $123/22$ in this case. Therefore this is not a solution 
mod $2$. Indeed, no solution mod $2$ exists, since the denominator of the $x$ coordinate is divisible by $2$. For the $y$ coordinate we obtain
\begin{verbatim}  
  -- y coordinate
  recoverQQ(-526048,2*13*31*37*41)
  
          77
  o22 = - --
          43
\end{verbatim}
As a last step we substitute this $\QQ$-point into the original equations.
\begin{verbatim}
  sub(I, matrix{{123/22,-77/43}})
  o24 = ideal (0, 0)
\end{verbatim}
This shows that we have indeed found a solution over $\QQ$. Notice also that
as argued above one of the denominators is divisible by $2$ and the other by $11$.
\end{experiment}

\begin{rem} $\quad$

\begin{enumerate}
\item The assumption that we have a unique solution over $\QQ$ is not as restrictive as it might seem. If we have for example $2$  solutions, then at least the line through them is unique. More generally,
if the solution set over $\QQ$ lies on $k$ polynomials of degree $d$ then the corresponding
point in the Grassmannian $\GG(k,{ d+n \choose n})$ is unique.
\item Even if we do not have isolated solutions, we can use this method to find the polynomials
of $\rad(I(X))$, at least if the polynomials are of small degree. 
\item For this method we do not need explicit equations, rather an algorithm that decides whether a point lies on $X$ is enough. This is indeed an important distinction. It is for example easy to check whether a given hypersurface is singular, but very difficult to give an explicit discriminant polynomial in the coefficients of the hypersurface that vanishes if and only if it is singular. 
\end{enumerate}
\end{rem}

Before we finish this tutorial by looking at a very nice application of this method by Oliver Labs, 
we will look briefly at a second situation in which we can find explicit solutions over $\QQ$. I learned this method from Noam Elkies in his talk at the Clay Mathematics Institute Summer School ``Arithmetic geometry" 2006.

\begin{figure}
\includegraphics*[width=10cm]{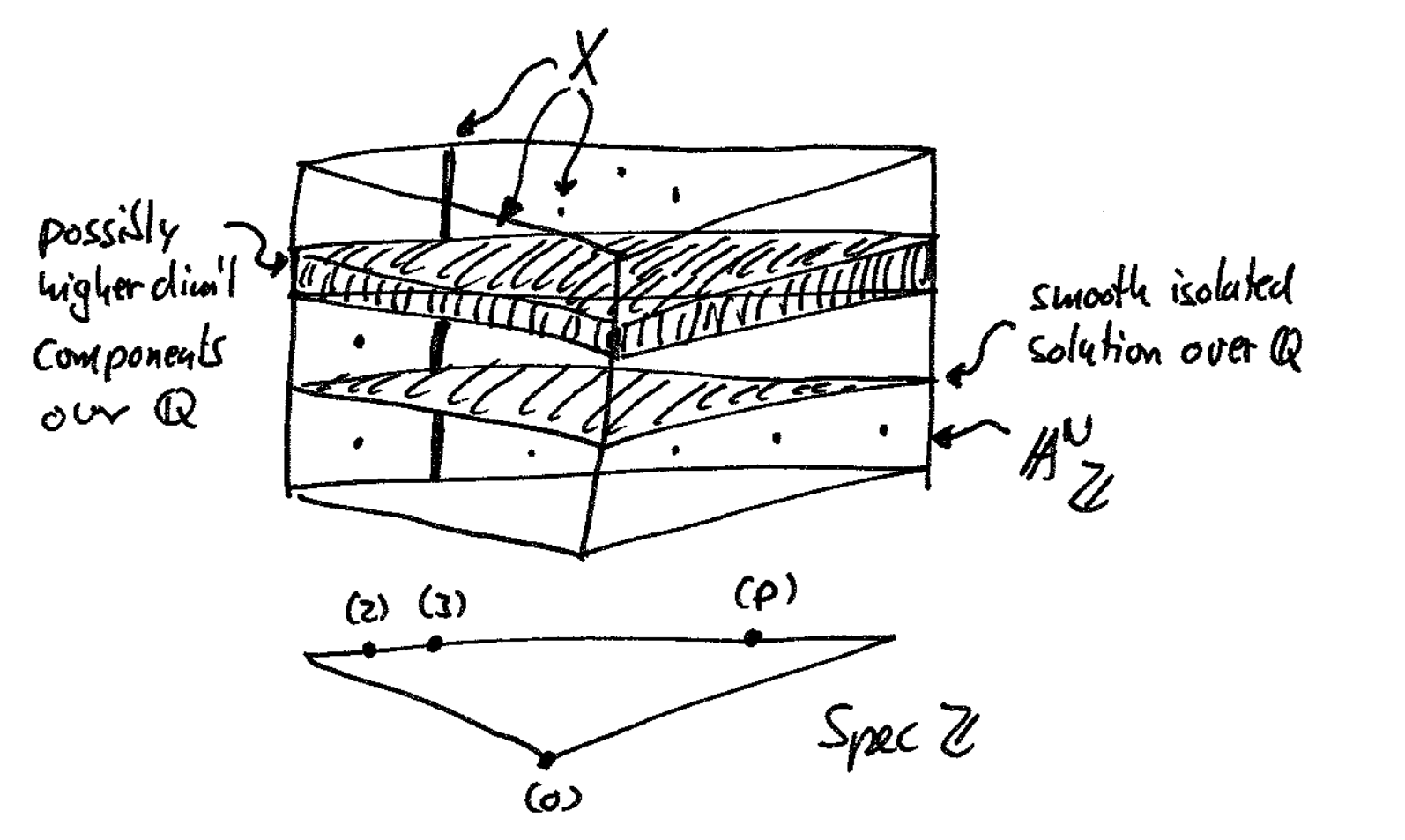}
\caption{A scheme $X$ over $\spec \ZZ$ with a smooth isolated solution over $\QQ$.}
\label{fLiftElkies}
\end{figure}

\begin{alg}
Assume that $X$ has a {\sl smooth} point $x$ over $\QQ$ that is isolated over $\bar{\QQ}$ 
as depicted in Figure \ref{fLiftElkies}, and that $p$ is a prime that does not
divide the denominators of the coordinates of $x$. Then we can find this point as follows:
\begin{enumerate}
\item Reduce mod $p$ and test all points.
\item Calculate the tangent spaces at the found points. If the dimension of such a tangent space
is $0$ then the corresponding point is smooth and isolated.
\item Lift the point mod $p^k$ with $k$ large using $p$-adic Newton iteration, as explained in Prop 
\ref{pNewton}.
\end{enumerate}
\end{alg}
 
 \begin{prop} \label{pNewton}
 Let $a \in \AZ^n_{\ZZ}$ be a solution of $$f_1(a)=\dots=f_n(a) = 0 \mod p^k$$ and 
 assume that the Jacobian matrix $J = \left(\frac{d f_i}{d x_j}\right)$ is invertible at $a$ mod $p$. Then
 \[
 	a' = a - (f_1(a),\dots,f_n(a)) J(a)^{-1}
\]
is a solution mod $p^{2k}$.
\end{prop}

\begin{proof}
Use the Taylor expansion as in the proof of Newton iteration.
\end{proof}

\begin{experiment}
Let's solve the equations of Experiment \ref{eLiftZZ} using $p$-adic Newton iteration. For this we
need some functions for modular calculations:
\begin{verbatim}
  -- calculate reduction of a matrix M mod n
  modn = (M,n) -> (
       matrix apply(rank target M, i-> 
       	  apply(rank source M,j-> M_j_i-round(M_j_i/n)*n)))

  -- divide a matrix of integers by an integer
  -- (in our application this division will not have a remainder)
  divn = (M,n) -> (
       matrix apply(rank target M, i-> 
       	  apply(rank source M,j-> M_j_i//n)))

  -- invert number mod n
  invn = (i,n) -> (
       c := gcdCoefficients(i,n);
       if c#0 == 1 then c#1 else "error"
       )

  -- invert a matrix mod n
  -- M a square matrix over ZZ
  -- (if M is not invertible mod n, then 0 is returned)
  invMatn = (M,n) -> (
       Mn := modn(M,n);
       MQQ := sub(Mn,QQ);
       detM = sub(det Mn,QQ);
       modn(invn(sub(detM,ZZ),n)*sub(detM*MQQ^-1,ZZ),n)
       )
\end{verbatim}
With this we can implement Newton iteration. We will represent  a point  by a pair $(P,eps)$ with
$P$ a matrix of integers that is a solution modulo $eps$.
\begin{verbatim}
  -- (P,eps)  an approximation mod eps (contains integers)
  -- M        affine polynomials (over ZZ)
  -- J        Jacobian matrix (over ZZ)
  -- returns an approximation (P,eps^2)
  newtonStep = (Peps,M,J) -> (
       P := Peps#0;
       eps := Peps#1;
       JPinv := invMatn(sub(J,P),eps);
       correction := eps*modn(divn(sub(M,P)*JPinv,eps),eps);
       {modn(P-correction,eps^2),eps^2}
       )

  -- returns an approximation mod Peps^(2^num)
  newton = (Peps,M,J,num) -> (
    i := 0;
    localPeps := Peps;
    while i < num do (
      localPeps = newtonStep(localPeps,M,J);
      print(localPeps);
      i = i+1;
    );
    localPeps
  )
\end{verbatim}
We now consider equations of Example \ref{eLiftZZ}
\begin{verbatim}
  I = ideal (-8*x^2-x*y-7*y^2+5238*x-11582*y-7696,
     4*x*y-10*y^2-2313*x-16372*y-6462)   
\end{verbatim}
their Jacobian matrix
\begin{verbatim}
  J = jacobian(I)    
\end{verbatim}
and their solutions over $\FF_7$:
\begin{verbatim}
  apply(allPoints(I,7),Pseq -> (
	    P := matrix {toList Pseq};
	    (P,0!=det modn(sub(J,P),7))
	  ))

  o25 = {(| 2 3 |, true), (| 5 5 |, true)}
\end{verbatim}
Both points are isolated and smooth over $\FF_7$ so we can apply $p$-adic
Newton iteration to them. The first one lifts to the solution found in Experiment \ref{eLiftZZ}:
\begin{verbatim}
  newton((matrix{{2,3}},7),gens I, J,4)

  {| 9 10 |, 49}
  {| -1167 -774 |, 2401}
  {| 1234 -774 |, 5764801}
  {| 1234 -774 |, 33232930569601}
\end{verbatim}
while the second point probably does not lift to $\ZZ$:
\begin{verbatim}
  newton((matrix{{5,5}},7),gens I, J,4)
  
  {| 5 -9 |, 49}
  {| -926 334 |, 2401}
  {| 359224 -66894 |, 5764801}
  {| 11082657337694 -9795607574104 |, 33232930569601}
\end{verbatim}
\end{experiment}

\begin{rem}
Noam Elkies has used this method to find interesting elliptic fibrations over $\QQ$. 
See for example \cite[Section III, p. 11]{ElkiesOberwolfach}.
\end{rem}

\begin{rem}
The Newton method is much faster than lifting by Chinese remaindering, since we only need to find
one smooth point in one characteristic. Unfortunately, it does not work if we cannot calculate tangent spaces.  An application where this happens
is discussed in the next section.
\end{rem}

%%%%%%%%%%%%%%%%%%%%%%%%%%%%%%%%%%%%%%%%%%%%%
\section{Surfaces with Many Real Nodes} \label{sNodes}
%%%%%%%%%%%%%%%%%%%%%%%%%%%%%%%%%%%%%%%%%%%%%
\nosubsections
A very nice application of finite field experiments with beautiful characteristic zero results was
done by Oliver Labs in his thesis \cite{labsPhD}. We look at his ideas and results in this section.

\newcommand{\pthreeR}{\PP^3_{\RR}}
Consider an algebraic surface $X  \subset \pthreeR$ of degree $d$ and denote by $N(X)$ the
number of real nodes of $X$. A classical question of real algebraic geometry is to determine
the maximal number of nodes a surface of degree $d$ can have. We denote this number by
\[
	\mu(d) := \max\{N(X) \suchthat X \subset \pthreeR \wedge \deg X = d\}.
\]
Moreover one would like to find explicit equations for surfaces $X$ that do have
$\mu(d)$ real nodes. The cases $\mu(1) = 0$ and $\mu(2)=1$, i.e the plane and the quadric cone,
have been known since antiquity. 

\begin{figure}
\includegraphics*[width=6cm]{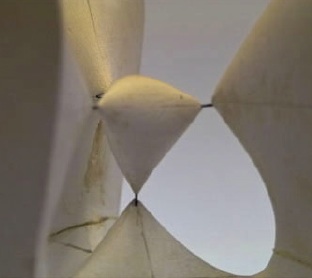}
\caption{Historic plaster model of the Cayley Cubic as displayed in the mathematical instute of the university of G\"ottingen}
\label{f136}
\end{figure}

\begin{figure}
\includegraphics*[width=6cm]{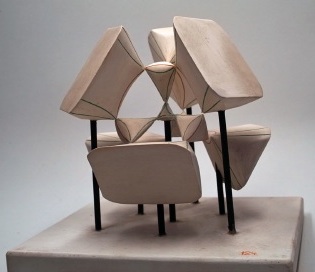}
\caption{A Kummer surface with 16 nodes.}
\label{f124}
\end{figure}

Cayley \cite{CayleyCubic} and Sch\"afli \cite{SchlaefliClassification} solved
$\mu(3) = 4$, while Kummer proved $\mu(4)=16$ in \cite{KummerClassification}. Plaster models of a Cayley-Cubic and a Kummer-Quartic
are on display in the G\"ottingen Mathematical Institute as numbers $124$ and $136$, see Figure \ref{f136} and \ref{f124}. These pictures many other are available at 

\begin{center}
\href{http://www.uni-math.gwdg.de/modellsammlung/}{http://www.uni-math.gwdg.de/modellsammlung.}
\end{center}

For the case $d=5$, Togliatti proved
in \cite{TogliattiQuintic} that quintic surfaces with $31$ nodes exist. One such surface is depicted in Figure \ref{fTogliatti}.
It took  $40$ years before Beauville \cite{Beauville31} finally proved that $31$ is indeed the maximal possible number. 

\begin{figure}
\includegraphics*[width=6cm]{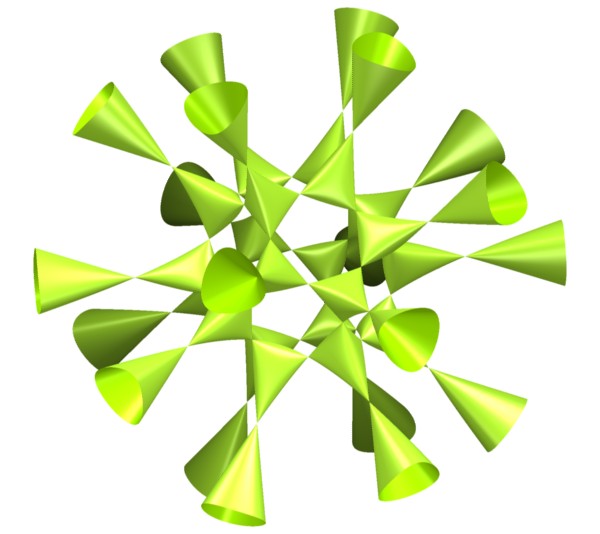}
\caption{The Barth sextic}
\label{fBarth}
\end{figure}

In 1994 Barth \cite{BarthSextic} found the beautiful sextic with the icosahedral symmetry and $65$ nodes
shown in Figure \ref{fBarth}. Jaffe and Rubermann proved in \cite{not66} that no sextics with $66$ or more nodes exist. 

For $d=7$ the problem is still open. By works of Chmutov \cite{ChmutovBound},  Breske/Labs/van Straten \cite{realArrNodes} and Varichenko \cite{VarchenkoBound} we only know
$93 \le \mu(7) \le 104$. For large $d$ Chmutov and Breske/Labs/van Straten show
$$\mu(d) \ge \frac{5}{12} d^3 +  \text{lower order terms},$$
while Miyaoka \cite{MiyaokaBound} proves
$$\mu(d) \le \frac{4}{9} d^3 + \text{lower order terms}.$$

Here we explain how Oliver Labs found a new septic with many nodes, using
finite field experiments \cite{LabsSeptic}. 

\begin{experiment}
The most naive approach to find septics with many nodes is to look at random surfaces
of degree $7$ in some small characteristic:
\begin{verbatim}
  -- Calculate milnor number for hypersurfaces in IP^3
  -- (for nonisolated singularities and smooth surfaces 0 is returned)
  mu = (f) -> (
       J := (ideal jacobian ideal f)+ideal f;     
       if 3==codim J then degree J else 0
       )

  K = ZZ/5        -- work in char 5
  R = K[x,y,z,w]  -- coordinate ring of IP^3

  -- look at 100 random surfaces
  time tally apply(100, i-> mu(random(7,R)))     
\end{verbatim}
After about $18$ seconds we find
\begin{verbatim}
  o4 = Tally{0 => 69}
           1 => 24
           2 => 5
           3 => 1
           4 => 1
\end{verbatim}
which is still far from $93$ nodes. Since having an extra node is a codimension-one condition, a rough estimation gives that  we would have to search $5^{89}  \approx 1.6 \times 10^{62}$ times longer to find $89$ more nodes in characteristic $5$. 
\end{experiment}

\begin{figure}
\includegraphics*[width=6cm]{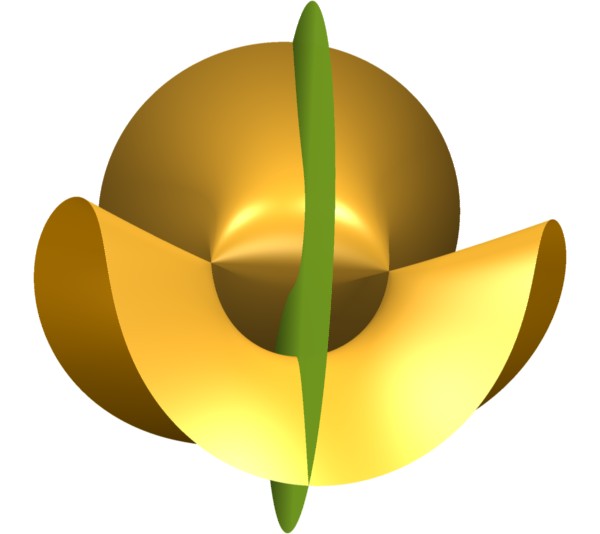}
\caption{A mirror symmetric cubic}
\label{fMirrorsym}
\end{figure}

One classical idea to find surfaces with many nodes, is to use symmetry. If for example
 we only look at mirror symmetric surfaces, we obtain singularities in pairs, as depicted in Figure \ref{fMirrorsym}.

 \begin{experiment}
 We look at $100$ random surfaces that are symmetric with respect to the $x=0$ plane
 \begin{verbatim}
  -- make a random f mirror symmetric
  sym = (f) -> f+sub(f,{x=>-x})

  time tally apply(100, i-> mu(sym(random(7,R))))   

  o6 = Tally{0 => 57}
           1 => 10
           2 => 11
           3 => 9
           4 => 4
           5 => 3
           6 => 3
           7 => 1
           9 => 1
           13 => 1
\end{verbatim}
Indeed, we obtain more singularities, but not nearly enough.
\end{experiment}

The symmetry approach works best if we have a large symmetry group. In the $d=7$ case
Oliver Labs used the $D_7$ symmetry of the $7$-gon. If $D_7$ acts on $\PP^3$ with symmetry axis $x=y=0$ one can use representation theory to find
a $7$-dimensional family of $D_7$-invariant $7$-tic we use in the next experiment.

\begin{experiment} \label{e7gon3D}
Start by considering the cone over a $7$-gon given by
$$
P = 2^6 \prod_{j=0^6} \left( \cos\left(\frac{2\pi j}{7} \right)x +  \sin\left(\frac{2\pi j}{7}\right) y-z\right),
$$
which can be expanded to 
\begin{verbatim} 
  P = x*(x^6-3*7*x^4*y^2+5*7*x^2*y^4-7*y^6)+
       7*z*((x^2+y^2)^3-2^3*z^2*(x^2+y^2)^2+2^4*z^4*(x^2+y^2))-
       2^6*z^7
\end{verbatim}
Now parameterize $D_7$ invariant septics $U$ that contain a double cubic.
\begin{verbatim}
  S = K[a1,a2,a3,a4,a5,a6,a7]
  RS = R**S  -- tensor product of rings
  U = (z+a5*w)*
      (a1*z^3+a2*z^2*w+a3*z*w^2+a4*w^3+(a6*z+a7*w)*(x^2+y^2))^2
\end{verbatim}
We will look at random sums of the form $P+U$ using
\begin{verbatim}
  randomInv = () -> (
       P-sub(U,vars R|random(R^{0},R^{7:0}))
       )
\end{verbatim}
Let's try $100$ of these
\begin{verbatim}
  time tally apply(100, i-> mu(randomInv()))   

o9 = Tally{63 => 48}
           64 => 6
           65 => 4
             ...
           136 => 1
           140 => 1
\end{verbatim}

\begin{figure}
\includegraphics*[width=6cm]{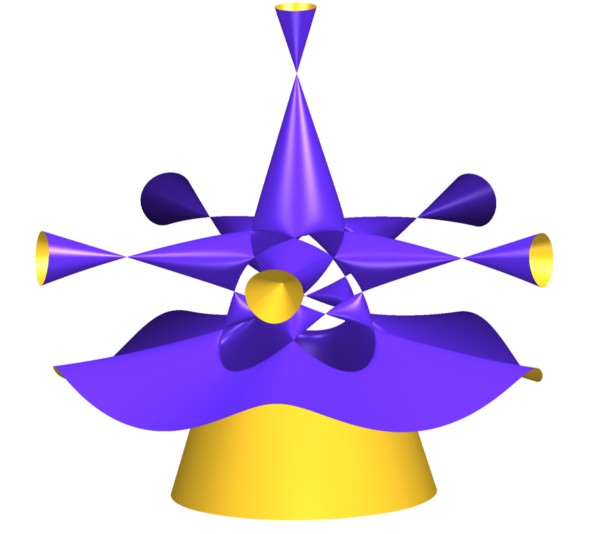}
\caption{A Togliatti quintic}
\label{fTogliatti}
\end{figure}

Unfortunately, this looks better than it is, since many of the surfaces with high Milnor
numbers have singularities that are not ordinary nodes. We can detect this by looking at the Hessian matrix which has rank $\ge 3$ only at smooth points and ordinary nodes. The following
function returns the number of nodes of $X =V(f)$ if all nodes are ordinary and $0$ otherwise.
\begin{verbatim}
  numA1 = (f) -> (
     -- singularities of f
     singf := (ideal jacobian ideal f)+ideal f;     
     if 3==codim singf then (
	       -- calculate Hessian
     	 Hess := diff(transpose vars R,diff(vars R,f));
	       ssf := singf + minors(3,Hess); 
	       if 4==codim ssf then degree singf else 0
	       )
	     else 0
   )
\end{verbatim}
With this we test another $100$ examples:
\begin{verbatim}
  time tally apply(100, i-> numA1(randomInv()))   

  o12 = Tally{0 => 28 }
            63 => 51
            64 => 13
            65 => 1
            70 => 6
            72 => 1
\end{verbatim}
which takes about $30$ seconds. Notice that most surfaces have $N(X)$ a multiple of $7$ as expected from the symmetry. 
\end{experiment}

\begin{figure}
\includegraphics*[width=6cm]{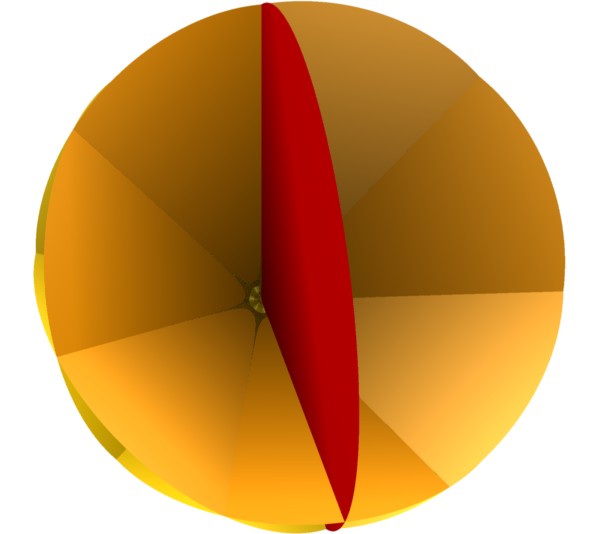}
\caption{Intersection of the $7$-gon with a perpendicular hypersurface}
\label{f7gon}
\end{figure}

To speed up these calculations Oliver Labs intersects the surfaces $X=V(P+U)$ with the hyperplane $y=0$ see Figure \ref{f7gon}. Since the operation of $D_7$ moves this hyperplane to $7$ different positions, every singularity of the intersection curve $C$ that does not lie on the symmetry axis corresponds to $7$ singularities  of $X$. Singular points on $C$ that do lie on the symmetry axis contribute only one node to the singularities of $X$. Using the symmetry of the construction one can show that for surfaces $X$ with only ordinary double points all singularities are obtained this way \cite[p.\,18, Cor.\,2.3.10]{EndrassPHD}, \cite[Lemma 1]{LabsSeptic}.

\begin{experiment}
We now look at $10000$ random $D_7$-invariant surfaces and their intersection curves with $y=0$.  We estimate the number of nodes on $X$ from the number of nodes on $C$ and return the point in the parameter space of $U$ if this number is large enough.
\begin{verbatim}
use R
time tally apply(10000,i-> (
	  r := random(R^{0},R^{7:0});
	  f := sub(P-sub(U,vars R|r),y_R=>0);
	  singf := ideal f + ideal jacobian ideal f;
	  if 2 == codim singf then (
	       -- calculate Hessian
     	 Hess := diff(transpose vars R,diff(vars R,f));
	       ssf := singf + minors(2,Hess); 
	       if 3==codim ssf then (
		          d := degree singf;
		          -- points on the line x=0
	          singfx := singf+ideal(x);
	 	        dx := degree singfx;
	          if 2!=codim singfx then dx=0;
	 	        d3 = (d-dx)*7+dx;		   
	          (d,d-dx,dx,d3,if d3>=93 then r)
	       ) 
	       else -1
	     )
	  ))
\end{verbatim}
In this way we find
\begin{verbatim}
o16 = Tally{(9, 9, 0, 63, ) => 5228                    
            (10, 9, 1, 64, ) => 731
             .....
            (16, 14, 2, 100, | 1 2 2 1 1 0 1 |) => 1
            -1 => 3071
            null => 8
\end{verbatim}
It remains to check whether the found $U$ really gives rise to surfaces with $100$ nodes
\begin{verbatim}
  f = P-sub(U,vars R|sub(matrix{{1,2,2,1,1,0,1}},R))
  numA1(f)
  
  o18 = 100
\end{verbatim}
This proves that there exists a surface with $100$ nodes over $\FF_5$.
\end{experiment}

\begin{figure}
\includegraphics*[width=6cm]{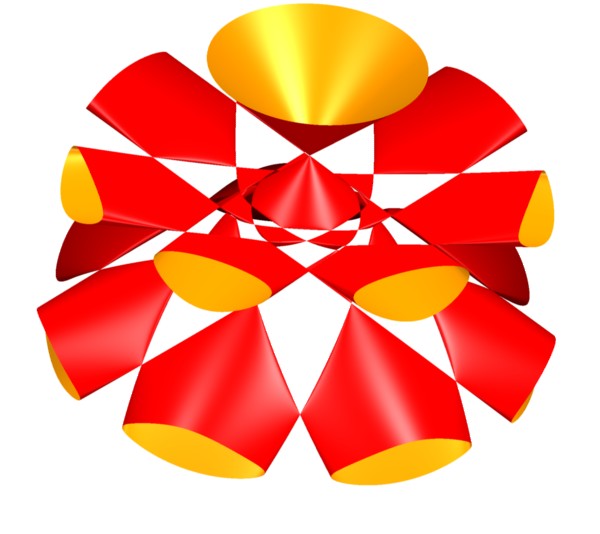}
\caption{The Labs septic}
\label{fLabs}
\end{figure}

Looking at other fields one finds that $\FF_5$ is a special case. In general one only finds
surfaces with $99$ nodes. To lift these examples to characteristic zero, Oliver Labs analyzed the
geometry of the intersection curves of the $99$-nodal examples and found that 
\begin{enumerate}
\item All such intersection curves decompose into a line and a 6-tic.
\item The singularities of the intersection curves are in a special position that can be explicitly described (see \cite{LabsSeptic} for details)
\end{enumerate}
These geometric properties imply (after some elimination) that there exists an $\alpha$ such that
\begin{align*}
\alpha_1 &= \alpha^7+7\alpha^5-\alpha^4+7\alpha^3-2\alpha^2-7\alpha-1 \\
\alpha_2 &= (\alpha^2+1)(3\alpha^5+14\alpha^4-3\alpha^2+7\alpha-2) \\
\alpha_3 &= (\alpha^1+1)^2(3\alpha^3+7\alpha-3) \\
\alpha_4&= (\alpha(1+\alpha^2)-1)(1+\alpha^2)^2 \\
\alpha_5 &= -\frac{\alpha^2}{1+\alpha^2} \\
\alpha_6 &= \alpha_7 = 1
\end{align*}
It remained to determine which $\alpha$ lead to $99$-nodal septics. 
Experiments over many primes show that there are at most $3$ such $\alpha$. Over primes with exactly $3$ solutions, Oliver Labs represented
them as zeros of a degree $3$ polynomial. By using the Chinese remaindering method, he lifted the
coefficients of this polynomial to characteristic $0$ and obtained 
\[
	7\alpha^3+7\alpha+1=0.
\] 
This polynomial has exactly one real solution, and with this $\alpha$ one can calculate this time
over $\QQ(\alpha)$ that the resulting septic has indeed $99$ real nodes. Figure \ref{fLabs} shows the inner part of this surface.

A movie of this and many other surfaces in this section can by found on my home page 
\begin{center}
\href{www.iag.uni-hannover.de/~bothmer/goettingen.php}{www.iag.uni-hannover.de/$\tilde{\,\,\,}$bothmer/goettingen.php},
\end{center}
on the home page of Oliver Labs
\begin{center}
\href{http://www.algebraicsurface.net/}{http://www.algebraicsurface.net/},
\end{center}
or on youTube.com 
\begin{center}
\href{http://www.youtube.com/profile?user=bothmer}{http://www.youtube.com/profile?user=bothmer}.
\end{center}
The movies and the surfaces in this article were produced using the 
public domain programs \verb=surf= by Stefan Endra\ss \, \cite{surf} and \verb=surfex= by Oliver Labs \cite{surfex}.

\begin{appendix}

\section{Selected Macaulay Commands} \label{Amacaulay}

Here we review some \verb!Macaulay 2! commands used in this tutorial. Lines starting with ``i" are input lines, while lines starting with ``o" are output lines. For more detailed explanations 
we refer to the online help of \verb!Macaulay2! \cite{M2}.

\subsection{apply}

This command applies a function to a list. In Macaulay 2 this is often used to generate loops.
\begin{verbatim}
  i1 : apply({1,2,3,4},i->i^2)
  o1 = {1, 4, 9, 16}
  o1 : List
\end{verbatim}
The list $\{0,1,\dots,n-1\}$ can be abbreviated by $n$:
\begin{verbatim}
  i2 : apply(4,i->i^2)
  o2 = {0, 1, 4, 9}
  o2 : List
\end{verbatim}

\subsection{map}

With \verb|map(R,S,m)| a map from $S$ to $R$ is produced. The matrix $m$ over $S$ contains the images of the variables of $R$: 
\begin{verbatim}
  i1 : f = map(ZZ,ZZ[x,y],matrix{{2,3}})
  o1 = map(ZZ,ZZ[x,y],{2, 3})
  o1 : RingMap ZZ <--- ZZ[x,y]

  i2 : f(x+y)
  o2 = 5
\end{verbatim}
If no matrix is given, all variables to variables of the same name or to zero.
\begin{verbatim}
  i3 : g = map(ZZ[x],ZZ[x,y])
  o3 = map(ZZ[x],ZZ[x,y],{x, 0})
  o3 : RingMap ZZ[x] <--- ZZ[x,y]

  i4 : g(x+y+1)
  o4 = x + 1
  o4 : ZZ[x]
\end{verbatim}

\subsection{random}

This command can be used either to construct random matrices
\begin{verbatim}
  i1 : K = ZZ/3

  o1 = K
  o1 : QuotientRing

  i2 : random(K^2,K^3)

  o2 = | 1 0  -1 |
       | 1 -1 1  |
               2       3
  o2 : Matrix K  <--- K
\end{verbatim}
or to construct random homogeneous polynomials of given degree
\begin{verbatim}
  i3 : R = K[x,y]
  
  o3 = R
  o3 : PolynomialRing

  i4 : random(2,R)
  
        2          2
  o4 = x  + x*y - y
  o4 : R
\end{verbatim}

\subsection{sub}

This command is used to substitute values for the variables of a ring:
\begin{verbatim}
  i1 : K = ZZ/3
 
  o1 = K
  o1 : QuotientRing
 
  i2 : R = K[x,y]
 
  o2 = R
  o2 : PolynomialRing

  i3 : f = x*y
  
  o3 = x*y
  o3 : R

  i4 : sub(f,matrix{{2,3}})
  
  o4 = 6
\end{verbatim}
Another application is the transfer a polynomial, ideal or matrix from one ring $R$ to another ring $S$ that has some variables in common with $R$
\begin{verbatim}
  i5 : S = K[x,y,z]
 
  o5 = S
  o5 : PolynomialRing
  
  i6 : sub(f,S)
  
  o6 = x*y
  o6 : S
\end{verbatim}

\subsection{syz}

The command is used here to calculate a presentation for the kernel of a 
matrix:

\begin{verbatim}
  i1 : M = matrix{{1,2,3},{4,5,6}}
  
  o1 = | 1 2 3 |
       | 4 5 6 |
                2        3
  o1 : Matrix ZZ  <--- ZZ

  i2 : syz M
  
  o2 = | -1 |
       | 2  |
       | -1 |
                3        1
  o2 : Matrix ZZ  <--- ZZ
\end{verbatim}

\subsection{tally}

With \verb|tally| one can count how often an element appears in a list:
\begin{verbatim}
  i1 : tally{1,2,1,3,2,2,17}
  
  o1 = Tally{1 => 2 }
             2 => 3
             3 => 1
             17 => 1
  o1 : Tally
\end{verbatim}

\section{Magma Scripts (by Stefan Wiedmann)} \label{Amagma}

Stefan Wiedmann  \cite{natoweb} has translated the Macaulay 2 scripts of this article to Magma. 
Here they are:

\begin{expB11} Evaluate a given polynomial in 700 random points.
\begin{verbatim}
K := FiniteField(7);               //work over F_7
R<x,y,z,w> := PolynomialRing(K,4); //Polynomialring in 4 variables over F_7
K4:=CartesianPower(K,4);           //K^4
F := x^23+1248*y*z*w+129269698;    //a polynomial

M := [Random(K4): i in [1..700]];  //random points
T := {*Evaluate(F,s): s in M*};

Multiplicity(T,0);                 //Results with muliplicity
\end{verbatim}
\end{expB11}

\begin{expB12} Evaluate a product of two polynomials in 700 random points
\begin{verbatim}
K := FiniteField(7);                //work over F_7
R<x,y,z,w> := PolynomialRing(K,4);  //AA^4 over F_7
K4:=CartesianPower(K,4);            //K^4

F := x^23+1248*y*z*w+129269698;     //a polynomial
G := x*y*z*w+z^25-938493+x-z*w;     //a second polynomial
H := F*G;

M := [Random(K4): i in [1..700]];   //random points
T := {*Evaluate(H,s): s in M*};
T;
Multiplicity(T,0);
\end{verbatim}
\end{expB12}

\begin{expB114} Count singular quadrics.
\begin{verbatim}
K := FiniteField(7);
R<X,Y,Z,W> := PolynomialRing(K,4); 
{* Dimension(JacobianIdeal(Random(2,R,0))) : i in [1..700]*};
\end{verbatim}
\end{expB114}

\begin{expB118} Count quadrics with $\dim > 0$ singular locus
\begin{verbatim}
function findk(n,p,k,c)
    //Search until k singular examples of codim at most c are found,
    //p prime number, n dimension
    K :=  FiniteField(p);
    R := PolynomialRing(K,n);
    trials := 0;
    found := 0;
    while found lt k do
      Q := Ideal([Random(2,R,0)]);
      if c ge n - Dimension(Q+JacobianIdeal(Basis(Q))) then
        found := found + 1;
      else
        trials := trials + 1;
      end if;
    end while;
    print "Trails:",trials;
    return trials;
end function;

k := 50;

time L1 := [[p,findk(4,p,k,2)] : p in [5,7,11]];
L1;

time findk(4,5,50,2);
time findk(4,7,50,2);
time findk(4,11,50,2);

function slope(L)
    //calculate slope of regression line by
    //formula form [2] p. 800
    xbar := &+[L[i][1] : i in [1..#L]]/#L;
    ybar := &+[L[i][2] : i in [1..#L]]/#L;
    return &+[(L[i][1]-xbar)*(L[i][2]-ybar): i in [1..3]]/
      &+[(L[i][1]-xbar)^2 : i in [1..3]];
end function;

//slope for dim 1 singularities
slope([[Log(1/x[1]), Log(k/x[2])] : x in L1]);
\end{verbatim}
\end{expB118}

\begin{expB21} Count points on a reducible variety.
\begin{verbatim}
K := FiniteField(7);
V := CartesianPower(K,6);
R<x1,x2,x3,x4,x5,x6> := PolynomialRing(K,6);

//random affine polynomial of degree d
randomAffine := func< d | &+[ Random(i,R,7) : i in [0..d]]>;

//some polynomials
F := randomAffine(2);
G := randomAffine(6);
H := randomAffine(7);

//generators of I(V(F) \cup V(H,G))
I := Ideal([F*G,F*H]);

//experiment
null := [0 : i in [1..#Basis(I)]];

t := {**};

for j in [1..700] do
    point := Random(V);
    Include(~t, null eq [Evaluate(Basis(I)[i],point) : i in [1..#Basis(I)]]);
end for;

//result
t;
\end{verbatim}
\end{expB21}

\begin{expB24} Count points and tangent spaces on a reducible variety.
\begin{verbatim}
K := FiniteField(7); //charakteristik 7
R<x1,x2,x3,x4,x5,x6> := PolynomialRing(K,6); //6 variables
V := CartesianPower(K,6);

//random affine polynomial of degree d
randomAffine := func< d | &+[ Random(i,R,7) : i in [0..d]]>;

//some polynomials
F := randomAffine(2);
G := randomAffine(6);
H := randomAffine(7);

//generators of I(V(F) \cup V(H,G))
I := Ideal([F*G,F*H]);

null := [0 : i in [1..#Basis(I)]];

//the Jacobi-Matrix
J := JacobianMatrix(Basis(I)); 
size := [NumberOfRows(J),NumberOfColumns(J)];

A := RMatrixSpace(R,size[1],size[2]);
B := KMatrixSpace(K,size[1],size[2]);

t := {**};

time
for j in [1..700] do
  point := Random(V);
  substitude := map< A -> B | x :-> [Evaluate(t,point): t in ElementToSequence(x)]>;
  if null eq [Evaluate(Basis(I)[i], point) : i in [1..#Basis(I)]] 
    then Include(~t, Rank(substitude(J)));
  else
    Include(~t,-1); 
  end if;
end for;

//result
t;
\end{verbatim}
\end{expB24}

\begin{expB29} $\quad$
\begin{verbatim}
K := FiniteField(7);                          //charakteristik 7
R<x1,x2,x3,x4,x5,x6> := PolynomialRing(K,6);  //6 variables
V := CartesianPower(K,6);

//consider an 5 x 5 matrix with degree 2 entries
r := 5; 
d := 2;
Mat := MatrixAlgebra(R,r);


//random matrix
M := Mat![Random(d,R,7) : i in [1..r^2]];

//calculate determinant and derivative w.r.t x1
time F := Determinant(M);
time F1 := Derivative(F,1);

//substitute a random point
point := Random(V);
time Evaluate(F1,point);

//calculate derivative with epsilon

Ke<e> := AffineAlgebra<K,e|e^2>; //a ring with e^2 = 0

Mate := MatrixAlgebra(Ke,r);

//the first unit vector
e1 := <>;
for i in [1..6] do
    if i eq 1 then
      Append(~e1,e);
    else
      Append(~e1,0);
    end if;
end for;

//point with direction
point1 := < point[i]+e1[i] : i in [1..6]>;

time Mate![Evaluate(x,point1): x in ElementToSequence(M)];
time Determinant(Mate![Evaluate(x,point1): x in ElementToSequence(M)]);

//determinant at 5000 random points

time
for i in [1..5000] do
    point := Random(V);                        //random point
    point1 := <point[i]+e1[i] : i in [1..6]>;  //tangent direction
    //calculate derivative
    _:=Determinant(Mate![Evaluate(x,point1): x in ElementToSequence(M)]);
end for;
\end{verbatim}
\end{expB29}

\begin{expB43} $\quad$
\begin{verbatim}
R<x,y> := PolynomialRing(IntegerRing(),2); //two variables

//the equations
F := -8*x^2-x*y-7*y^2+5238*x-11582*y-7696;
G := 4*x*y-10*y^2-2313*x-16372*y-6462;
I := Ideal([F,G]);

//now lets find the points over F_p
function allPoints(I,p)
    M := [];
    K := FiniteField(p);
    A := AffineSpace(K,2); 
    R := CoordinateRing(A);
    for pt in CartesianPower(K,2) do
        Ipt := Ideal([R|Evaluate(Basis(I)[k],pt) : k in [1..#Basis(I)]]);
        SIpt := Scheme(A,Ipt);
        if Codimension(SIpt) eq 0 then;
            Append(~M,pt);
        end if;
    end for;
    return M;
end function;

for p in PrimesUpTo(23) do
    print p, allPoints(I,p);
end for;

/*Chinese remaindering
given solutions mod m and n find
a solution mod m*n
sol1 = [n,solution]
sol2 = [m,solution]*/
function chinesePair(sol1,sol2)
    n := sol1[1];
    an := sol1[2];
    m := sol2[1];
    am := sol2[2];
    d,r,s :=   Xgcd(n,m);
    //returns d,r,s so that a*r + b*s is
    //the greatest common divisor d of a and b.
    amn := s*m*an+r*n*am;
    amn := amn - (Round(amn/(m*n)))*(m*n);
    if d eq 1 then
        return [m*n,amn];
    else
        print "m and n not coprime";
        return false;
    end if;
end function;

/*take a list {(n_1,s_1),...,(n_k,s_k)}
and return (n,a) such that
n=n_1* ... * n_k and
s_i = a mod n_i*/
function chineseList(L)
    //#L >= 2
    erg := L[1];
    for i in [2..#L] do
        erg := chinesePair(L[i],erg);
    end for;
    return erg;
end function;

//x coordinate   
chineseList([[2,0],[5,4],[17,10],[23,15]]);
//y coordinate
chineseList([[2,0],[5,1],[17,8],[23,8]]);

//test the solution
Evaluate(F,[1234,-774]);
Evaluate(G,[1234,-774]);
\end{verbatim}
\end{expB43}

\begin{expB46}  Rational recovery, as suggested in von zur Gathen in \cite[Section 5.10]{vzGathen}. Uses the functions \verb!allPoints! and \verb!chineseList! from Experiment B.4.3.
\begin{verbatim}
R<x,y> := PolynomialRing(IntegerRing(),2); //two variables

//equations
F := 176*x^2+148*x*y+301*y^2-742*x+896*y+768;
G := -25*x*y+430*y^2+33*x+1373*y+645;
I := Ideal([F,G]);

for p in PrimesUpTo(41) do
    print p, allPoints(I,p);
end for;

// x coordinate
chineseList([[2,1],[13,5],[31,7],[37,14],[41,0]]);
// y coordinate
chineseList([[2,0],[13,10],[31,22],[37,18],[41,23]]);

//test the solution
Evaluate(F,[138949,-526048]);
Evaluate(G,[138949,-526048]);

/*take (a,n) and calculate a solution to
r = as mod n
such that r,s < sqrt(n).
return (r/s)*/
function recoverQQ(a,n)
    r0:=a;
    s0:=1;
    t0:=0;
    r1:=n;
    s1:=0;
    t1:=1;
    r2:=0;
    s2:=0;
    t2:=0;
    k := Round(Sqrt(r1*1.0));
    while k le r1 do 
        q := r0 div r1;
        r2 := r0-q*r1;
        s2 := s0-q*s1;
        t2 := t0-q*t1;
        r0:=r1;
        s0:=s1;
        t0:=t1;
        r1:=r2;
        s1:=s2;
        t1:=t2;
    end while;
    return (r2/s2);
end function;

//x coordinate
recoverQQ(138949,2*13*31*37*41);
//y coordinate
recoverQQ(-526048,2*13*31*37*41);

//test the solution
Evaluate(F,[123/22,-77/43]);
Evaluate(G,[123/22,-77/43]);
\end{verbatim}
\end{expB46}

\begin{expB410} Lifting solutions using $p$-adic Newtoniteration (as suggested by N.Elkies). Uses the function \verb|allPoints| from Example B.4.3.
\begin{verbatim}
//calculate reduction of a matrix M mod n
function modn(M,n)
    return Matrix(Nrows(M),Ncols(M),[x - Round(x/n)*n : x in Eltseq(M)]);
end function;

//divide a matrix of integer by an integer
//(in our application this division will not have a remainder)
function divn(M,n)
    return Matrix(Nrows(M),Ncols(M),[x div n : x in Eltseq(M)]);
end function;

// invert number mod n
function invn(i,n)
    a,b := Xgcd(i,n);
    if a eq 1 then
        return b;
    else return false;
    end if;
end function;

//invert a matrix mod n
//M a square matrix over ZZ
//(if M is not invertible mod n, then 0 is returned)
function invMatn(M,n)
    Mn := modn(M,n);
    MQQ := MatrixAlgebra(RationalField(),Nrows(M))!Mn;
    detM := Determinant(Mn);
    if Type(invn(detM,n)) eq BoolElt then
        return 0;
    else
        return
            (MatrixAlgebra(IntegerRing(),Nrows(M))!
            (modn(invn(detM,n)*detM*MQQ^(-1),n)));
    end if;
end function;

//(P,eps) an approximation mod eps (contains integers)
//M       affine polynomials (over ZZ)
//J       Jacobian matrix (over ZZ)
//returns an approximation (P,eps^2)
function newtonStep(Peps,M,J)
    P := Peps[1];
    eps := Peps[2];
    JatP:=Matrix(Ncols(J),Nrows(J),[Evaluate(x,Eltseq(P)) : x in Eltseq(J)]);
    JPinv := invMatn(JatP,eps);
    MatP:= Matrix(1,#M,[Evaluate(x,Eltseq(P)) : x in Eltseq(M)]);
    correction := eps*modn(divn(MatP*Transpose(JPinv),eps),eps);
    return <modn(P-correction,eps^2),eps^2>;
end function;

//returns an approximation mod Peps^(2^num)
function newton(Peps,M,J,num)
    localPeps := Peps;
    for i in [1..num] do
	localPeps := newtonStep(localPeps,M,J);
	print localPeps;
    end for;
    return localPeps;
end function;

//c.f. example 4.3
R<x,y> := PolynomialRing(IntegerRing(),2); //two variables

//the equations
F := -8*x^2-x*y-7*y^2+5238*x-11582*y-7696;
G := 4*x*y-10*y^2-2313*x-16372*y-6462;
I := Ideal([F,G]);
J := JacobianMatrix(Basis(I));

Ap := allPoints(I,7);
for x in Ap do
    MatP := Matrix(Nrows(J),Ncols(J),[Evaluate(j,x): j in Eltseq(J)]);
    print x, (0 ne Determinant(MatP));
end for;

Peps := <Matrix(1,2,[2,3]),7>;
newton(Peps,Basis(I),J,4);

Peps := <Matrix(1,2,[5,5]),7>;
newton(Peps,Basis(I),J,4);
\end{verbatim}
\end{expB410}

\begin{expB51} Count singularities of random surfaces over $\FF_5$.
\begin{verbatim}
K := FiniteField(5);              //work in char 5
A := AffineSpace(K,4); 
R<x,y,z,w> :=CoordinateRing(A);   //coordinate ring of IP^3 

//Calculate milnor number
//(For nonisolated singularities and smooth surfaces 0 is returned)
function mu(f)
    SJ := Scheme(A,(Ideal([f])+JacobianIdeal(f)));
    if 3 eq Codimension(SJ) then
        return Degree(ProjectiveClosure(SJ));
    else
        return 0;
    end if;
end function;

//look at 100 random surfaces
M := {**};
time 
for i in [1..100] do
    f := Random(7,CoordinateRing(A),5);
    Include(~M,mu(f));
end for;

print "M:", M;
\end{verbatim}
\end{expB51}

\begin{expB52} Count singularities of mirror symmetic random surfaces. Uses the function $\verb|mu|$ from Example B.5.1.
\begin{verbatim}
K := FiniteField(5);              //work in char 5
A := AffineSpace(K,4); 
R<x,y,z,w> :=CoordinateRing(A);   //coordinate ring of IP^3 

//make a random f mirror symmetric
function mysym(f)
    return (f + Evaluate(f,x,-x));
end function;

//look at 100 random surfaces
M := {**};
time 
for i in [1..100] do
f:=R!mysym(Random(7,CoordinateRing(A),5));
Include(~M,mu(f));
end for;

print "M:", M;
\end{verbatim}
\end{expB52}

\begin{expB53} Count $A1$-singularities of $D_7$ invariant surfaces. Uses the function \verb|mu| from Experiment B.5.1.
\begin{verbatim}
K := FiniteField(5);              //work in char 5
A := AffineSpace(K,4); 
R<X,Y,Z,W> :=CoordinateRing(A);   //coordinate ring of IP^3 
RS<x,y,z,w,a1,a2,a3,a4,a5,a6,a7> := PolynomialRing(K,11);

//the 7-gon
P := X*(X^6-3*7*X^4*Y^2+5*7*X^2*Y^4-7*Y^6)
    +7*Z*((X^2+Y^2)^3-2^3*Z^2*(X^2+Y^2)^2+2^4*Z^4*(X^2+Y^2))-2^6*Z^7;

//parametrising invariant 7 tics with a double cubic
U := (z+a5*w)*(a1*z^3+a2*z^2*w+a3*z*w^2+a4*w^3+(a6*z+a7*w)*(x^2+y^2))^2;

//random invariant 7-tic
function randomInv()
    return (P - Evaluate(U,[X,Y,Z,W] cat [Random(K):i in [1..7]]));  
end function;

//test with 100 examples
M1 := {**};
time
for i in [1..100] do
Include(~M1,mu(randomInv()));
end for;

print "M1:", M1;

//singularities of f
function numA1(f)
    singf := Ideal([f])+JacobianIdeal(f);
    Ssingf := Scheme(A,singf);
    if 3 eq Codimension(Ssingf) then
        T := Scheme(A,Ideal([f]));
        //calculate Hessian
        Hess := HessianMatrix(T);
        ssf := singf + Ideal(Minors(Hess,3));
        if 4 eq Codimension(Scheme(A,ssf)) then
            return Degree(ProjectiveClosure(Ssingf));
        else
            return 0;
        end if;
    else
        return 0;
    end if;
end function;

//test with 100 examples
M2 := {**};
time
for i in [1..100] do
Include(~M2,numA1(randomInv()));
end for;

print "M2:", M2;

\end{verbatim}
\end{expB53}

\begin{expB54} Estimate number of $A1$-singularities by looking at $y=0$. Uses the function \verb|numA1| from Experiment B.5.3.
\begin{verbatim}
K := FiniteField(5); //work in char 5
A := AffineSpace(K,4); 
R<X,Y,Z,W> := CoordinateRing(A); //coordinate ring of IP^3 
RS<x,y,z,w,a1,a2,a3,a4,a5,a6,a7> := PolynomialRing(K,11);

//the 7-gon
P := X*(X^6-3*7*X^4*Y^2+5*7*X^2*Y^4-7*Y^6)
    +7*Z*((X^2+Y^2)^3-2^3*Z^2*(X^2+Y^2)^2+2^4*Z^4*(X^2+Y^2))-2^6*Z^7;

//parametrising invariant 7 tics with a double cubic
U := (z+a5*w)*(a1*z^3+a2*z^2*w+a3*z*w^2+a4*w^3+(a6*z+a7*w)*(x^2+y^2))^2;

//estimate number of nodes
function numberofsing()
    r := [Random(K): i in [1..7]];
    f := Evaluate(P - Evaluate(U,[X,Y,Z,W] cat r),Y,0);
    singf := Ideal([f])+ JacobianIdeal(f);
    Ssingf := Scheme(A,singf);
    if 2 eq Codimension(Ssingf) then
        //calculate Hessian
        S := Scheme(A,Ideal([f]));
        Hess := HessianMatrix(S);
        ssf := singf + Ideal(Minors(Hess,2));
        Sssf := Scheme(A,ssf);
        if 3 eq Codimension(Sssf) then
            d := Degree(ProjectiveClosure(Ssingf));
            //points on the line x=0
            singfx := singf + Ideal([X]);
            Ssingfx :=  Scheme(A,singfx);
            dx := Degree(ProjectiveClosure(Ssingfx));
            if 2 ne Codimension(Ssingfx) then
                dx := 0;
            end if;
            d3 := (d-dx)*7+dx;
            if d3 ge 93 then 
                return <d, d-dx, dx, d3>, r;
            else
                return <d, d-dx, dx, d3>, _;
            end if;
        else
            return <-1,0,0,0>,_;
        end if;
    else return <0,0,0,0>,_;
    end if;
end function;

M1 := {**};
M1hit := {**};
time
for i in [1..10000] do 
    a,b := numberofsing();
    if assigned(b) then
        Include(~M1hit,<a,b>);
    else
        Include(~M1,a);
    end if;
end for;

print "M1:", M1;
print "M1hit:", M1hit;

//test
f :=  P-Evaluate(U,[X,Y,Z,W,1,2,2,1,1,0,1]);
numA1(f);
\end{verbatim}
\end{expB54}

\end{appendix}

%\bibliographystyle{hunsrt} 
%\bibliography{../../../../../projekte/bib}

\def\cprime{$'$} \def\cprime{$'$}

 \end{document}